\magnification=1000
\baselineskip12.7pt
\hsize=13.7truecm
\vsize=22.5truecm
\leftskip=0truecm
\rightskip=0truecm
\hoffset=1,22truecm
\voffset=0truecm
\hfuzz=10pt
\pretolerance=500
\tolerance=1000
\brokenpenalty=5000

\topskip=2truecm

\def\D{\displaystyle}
\def\Ev{\vskip0,5cm}
\def\ev{\vskip0,25cm}

\def\min{\hbox{min}}
\def\max{\hbox{max}}

\def\n{\noindent}

\def\é{\hbox{\'e}}
\def\ê{\hbox{\^e}}
\def\è{\hbox{\`e}}
\def\à{\hbox{\`a}}
\def\ù{\hbox{\`u}}

\def\ti{\widetilde}
\def\g{\gamma}

\def\picture #1 by #2 (#3){\vbox to #2{\hrule width #1 height
0pt
\vfill\special{picture #3}}}
\def\scaledpicture #1 by #2 (#3 scaled
#4){{\dimen0=#1\dimen1=#2\divide\dimen0 by 1000
\multiply\dimen0 by #4\divide\dimen1 by 1000 \multiply\dimen1
by #4
\picture \dimen0 by \dimen1 (#3 scaled #4)}}

\font\bb=msbm10

\font\bbbb=msbm10 at 7pt
\def\d #1{\hbox{\bb #1}}

\def\r #1{\hbox{$\cal #1 $}}

\font\tm=cmr10
\font\tmm=cmr5

\font\Tm=cmr10 at 14pt
\font\tt=cmr9 at 8pt
\font\tcq=cmr10 at 5pt
\font\tcqq=cmr10 at 8pt

\font\tttsl=cmsl10 at 7pt

\font\LL=lasy10 at 14pt
\def\f{{\hbox{\hskip 5mm \LL 2}}}

\vglue0cm
\rm

\def\date{\the\day\  \ifcase\month\or janvier\or f\'evrier\or
mars\or avril\or mai \or juin \or juillet\or ao\^ut\or
septembre\or octobre\or novembre\or d\'ecembre\fi
\  {\the\year} }

\hbox to 13.5cm{\hrulefill}
\hskip-14mm\vtop{\hsize=13,5truecm  {\tttsl \baselineskip8pt
\date
\vfill} }\ \kern0cm\hfill\

\vbox{ \vskip1,5cm}
\centerline{\hskip0mm \hbox{\Tm ROTATION OF TRAJECTORIES OF LIPSCHITZ }}
\vskip3,5mm
\centerline{\hskip0mm \hbox{\Tm VECTOR FIELDS
       } }
\vskip10mm
\centerline{\hskip0mm Georges COMTE
\footnote{$(\dagger)$}{\tcqq Laboratoire J. A. Dieudonn\'e,
UMR CNRS 6621,
  Universit\'e de Nice Sophia-Antipolis,
28, avenue de Valrose, 06108 Nice Cedex 2, France.
\hskip1mm  e-mail: comte@math.unice.fr  },
Yosef YOMDIN\footnote{$(\ddagger)$}{
\tcqq Department of Mathematics, the Weizmann Institute of Science,
Rehovot 76100, Israel. \hskip1mm e-mail: yosef.yomdin@weizmann.ac.il.
Research supported by the ISF, Grant No. 264/02,
and by BSF, Grant No.2002243.} }
\vskip-0,5mm
 \centerline{ \hbox to 1.5cm{\hrulefill}}
\vskip3mm
\hskip-7mm \vtop{\hsize=13,6truecm \tolerance10000

{\tt \baselineskip=9pt
A{\tcq BSTRACT}. --  We prove that in finite time a trajectory of
a Lipschitz vector field in $\hbox{\bbbb R}^{\hbox{\tmm n}}$ can
not have infinite rotation around a given point.   This result
extends to the mutual rotation of two trajectories of a field in
$\hbox{\bbbb R}^{\hbox{\tmm 3}}$: this rotation is bounded from
above on any finite time interval. The bounds we give are only in
terms of the Lipschitz constant of the field and the length of the
time interval. \vfill} }\ \kern0cm\hfill\ \Ev

\vskip1,5cm

\centerline{\bf 1. Introduction}
\Ev

In this paper we investigate the tameness of a geometric behavior
of trajectories of vector fields: the rotation of such a
trajectory around a point or the mutual rotation of two
trajectories.

Of course a lot of results have been published on the geometry
of solutions of differential equations and trajectories of vector
fields, and we
simply cannot give an extensive review or even a bibliography
on the subject.
Nevertheless, we have at our disposal only few global theories
and general results concerning the tameness of trajectories.
In this very short introduction we just want to focus on two
of them.
\ev
In Gabrielov-Hovanskii's theory of Pfaffian sets, a Pfaffian function
 $f$ on an open set $\r U\subset \d R^n$,
is defined as a function that can be
written in the following way:
 $f(x)=P(x,c_1(x), \cdots, c_m(x))$, where $P$ is a polynomial and
 the $c_i$'s are analytic solutions of the polynomial
triangular differential system:
$$ Dc_i(x)=\D\sum_{j=1}^m P_{i,j}(x, c_1(x), \cdots,
c_i(x))\ dx_j, \ \
i\in \{1, \cdots, m\}. \eqno(*)$$
Then if we consider the Pfaffian structure, that is to say the smallest
structure containing the semi-Pfaffian sets,
it has been proved in [Wi], using
a Bezout type theorem of Hovanskii ([Ho 1], [Ho 2]), that this
structure is o-minimal (see [Co], [Dr], [Dr-Mi], [Sh]): for short
the number of connected components
of sets in such a structure is finite, and consequently no Pfaffian
curve (nor set) may infinitely oscillate or spiral around a point
and two
such curves have bounded mutual rotation.

This theory presents a large class of tame objects coming from
differential equations, but when we deal with vector fields, the
only way to a priori be sure that the trajectories belong to this
category is to assume that the trajectory is a Pfaffian function
itself, satisfying equation~$(*)$. We deduce from this assumption
that the field depends only on one variable. A too restrictive
hypothesis that obviously not allows a wild behavior for
trajectories.

On the other hand, starting with a given vector field and studying
the local geometry of its trajectories in a neighborhood of one of
its singularities, we know that this geometry is tame, following
[Ku-Mo] and [Ku-Mo-Pa], provided that the field is an analytic
gradient vector field: the rotation of a trajectory of an analytic
gradient vector field around one of its singularities is finite.
As a consequence the limit of the secant lines to the trajectory
exists (Thom's Gradient Conjecture). Let us notice that we do not
know wether a trajectory of an analytic gradient vector field lies
in some o-minimal structure, although we know that trajectories of
the gradient of a function definable in a given o-minimal
structure have finite length (see [Ku]). On the local behavior of
trajectories we would like to refer to the number of deep results
produced by the {\sl Spanish-Dijon School}, and specially the most
recent ones: [Ca-Mo-Ro],  [Ca-Mo-Sa 1,2,3], [Bl-Mo-Ro]. \ev

The aim of this  paper is to give  general results about rotation
of  trajectories of a vector field, with no restrictive
assumptions on the nature of the field, besides the Lipschitz
property which is a minimal hypothesis for the existence of
trajectories.

Of course, in this direction,  we cannot hope to treat
infinite time phenomena, as it is done in  [Ku-Mo-Pa] and the
Pfaffian theory, because
it is well-known that complete trajectories of polynomial vector fields in the
plane may spiral around a singularity of the field (see for
instance Remark of Section $3.1$). We thus obtain our bounds for
trajectories defined on a finite time interval.
For bounds obtained in the same spirit,
the reader may report to  [Gr-Yo], [Ho-Ya], [No-Ya 1,2], [Ya].
 \ev

The paper is organized as follows.

In Section $2$ we introduce
and compare
some notions of absolute and topological rotation for trajectories
 around
an affine subspace in $\d R^n$
or for two trajectories around each others in $\d R^3$.

In Section $3$ we first notice that the rotation of any
trajectory of a Lipschitz vector field around its stationary points
is bounded in terms of the Lipschitz constant (and the time
elapsed) only (Proposition $3.1$). The same is true for the rotation of any trajectory
around a linear invariant subspace of the field (Proposition $3.2$).
Moreover,  while the rotation velocity of a
trajectory around a {\it
non-stationary point} of the field may tend to infinity, we prove our
first main result: the ``total" rotation
around such a point is still bounded in terms of the Lipschitz
constant and the  time interval (Theorem $3.4$).
Our second main result is a consequent of the first one: we give
a uniform bound
for the rotation of any two trajectories of a given Lipschitz vector
field, in terms of the time interval and the Lipschitz constant
(Theorem $3.8$).
In contrast, we provide an easy example showing that the
rotation of a trajectory of $C^{\infty}$ vector field around a
non-invariant subspace may be infinite in finite time
(Example $3.3$).

\Ev
\centerline{\bf 2. Definition of signed and absolute rotation}
\Ev
{\bf 2.1. Rotation around an affine subspace}
\ev

First of all, we define an absolute rotation of the curve $\gamma$
in $\d R^n$ around the origin $0 \in \d R^n$.
Assuming that $\gamma$ does not pass through the origin, we can
define the spherical image of $\gamma$ as the curve $\sigma$ in the unit
sphere $S^{n-1}$, in the following way:

$$ \sigma(t)=\D {\gamma(t)\over \Vert\gamma(t)\Vert}. $$
\ev
{\bf Definition.}
 {\sl The absolute rotation } $R_{abs}(\gamma,0)$ of the curve
$\gamma$ around the origin is the length of $\sigma$, the spherical
image  of $\gamma$ in the unit sphere $S^{n-1}$.

\ev

We have the following lemma (easy to check) which
relates the rotation to the spherical part of the velocity of the
curve:
\ev
{\bf Lemma 2.1. --- }
{\sl
 Let
$\gamma^\prime(t)=\D{d\gamma(t)\over dt}$ be the velocity vector of
$\gamma:I\to \d R^n$, and let $\gamma^\prime_r(t)$ and $\gamma^\prime_s(t)$ be
the radial and the spherical components of this velocity vector.
Then the velocity of the spherical blowing-up $\sigma$ of $\gamma$ is:
$$ \sigma'(t)={\gamma'_s(t) \over \Vert \gamma(t)\Vert }.$$
As a consequence,
 the absolute rotation $R_{abs}(\gamma,0)$ is given by the
integral:
 $$ R_{abs}(\gamma,0)=\int_{t\in I} {{\Vert
\gamma^\prime_s(t)\Vert}\over {\Vert\gamma(t)\Vert}} \ dt.$$}
\ev
{\bf Remark.}
The absolute rotation $R_{abs}(\gamma,0)$
is invariant with respect to the monotone
reparametrizations of the curve $\gamma$ and in the same time
we want
our results  to be  about the geometry of the curve and not
depending of one of its parametrizations.
This is why in what
follows we implicitly consider  the geometric trajectory $\gamma(I)$
of the injective curve $\gamma:I\to \gamma(I)$ as the class of
$\gamma:I\to \d R^n$
modulo its injective reparametrizations.

\ev

For plane curves $\gamma$ in $\d R^2$ we can define their
``signed rotation" $R(\gamma,0)$ around the origin, which is, of
course, the usual rotation index. Indeed, in this case the unit
sphere  is a circle. Assuming that the orientation of
this circle (and of the curve $\gamma$) has been chosen, we can
define the signed rotation, essentially, by the same expression as
above:
\ev
{\bf Definition.}
 The {\sl signed rotation} $R(\gamma,0)$ of $\gamma$ around the origin
is defined by:
$$ R(\gamma,0)={1\over {2\pi}}\int_{\gamma}\pm {{ \Vert
\gamma^\prime_s(t)\Vert}\over {\Vert\gamma(t)\Vert}} dt. $$
Here the sign under the integral is chosen according to the direction
of the tangent to the unit circle vector $\gamma^\prime_s(t)$.

\ev
{\bf Remark.}
 Notice the normalization  by $2\pi$,
the unit  circle length, which appears here to make the linking
number, defined below, an integer. In  codimension greater or
equal to three we do not normalize the length of the spherical
curves.
\ev

An absolute rotation of the curve $\gamma$ in $\d R^n$
around a linear $k$-dimensional subspace ${\cal{L}} \subset
\d R^n$ is defined as follows: let ${\cal L}^\perp$ denote
the orthogonal subspace to $\cal L$. Let $\tilde \gamma$ be
the projection of  $\gamma$ on ${\cal L}^\perp$. Assuming
that $\gamma$ does not touch $\cal L$, we get $\tilde \gamma$ not
passing through the origin in ${\cal L}^\perp$.
\ev
{\bf Definition. } The {\sl absolute
rotation} $R_{abs}(\gamma,{\cal L})$ is defined as the absolute
rotation of the curve $\tilde \gamma$ in ${\cal L}^\perp$ around
the origin.
In the case of a linear subspace $\cal{L}$ of codimension $2$ in
$\d R^n$, a signed rotation $R(\gamma,{\cal L})$ of $\gamma$
around $\cal{L}$ is defined as the signed rotation of the curve
$\tilde \gamma$ in the plane ${\cal L}^\perp$ around the origin.

\ev

Of course the absolute rotation always bounds from above the
absolute value of the signed one.

\ev

For a closed curve $\gamma$ and for a subspace $\cal L$ of codimension $2$
the signed rotation $R(\gamma,{\cal L})$ is an integer, and it is
a topological invariant. For $\gamma$  non-closed   this signed
rotation $R(\gamma,{\cal L})$ may take non-integer values.
However, it is still an invariant of deformations of $\gamma$,
preserving the end points, and not touching $\cal{L}$.

\ev
{\bf 2.2. Rotation of two curves in $\d R^3$}
\ev
It is well known that the linking number of two closed curves in
$\d R^3$ can be defined via an integral expression,
the so-called Gauss integral (see, for example [Ar-Kh], [Du-Fo-No]).
This gives us a natural way to define also an absolute and a signed
rotation of two curves (closed or non-closed) one around the other.

For non-closed curves the
rotation defined in this (or any other) way cannot be metrically
or topologically invariant. But on the other
hand, the Gauss integral representation provides a powerful
analytic tool for its investigation.
 Our presentation in the next
subsection follows very closely the one given in [Du-Fo-No]. \ev
{\bf 2.2.1. The Linking Coefficient}
\ev
Consider a pair of smooth, closed, regular directed curves
 in $\d R^3$, which do not
intersect. We may assume them to be parametrized in
the following way:
 $\gamma_i:I_i\to \d R^3$, with $I_1, I_2$ two compact intervals.
 We denote our geometric curves by $\gamma_1$ and $\gamma_2$ (instead of
$\g_1(I_1), \g_2(I_2)$).

\ev
{\bf Definition.}
 The {\sl linking coefficient of the two curves} $\gamma_1 ,
\gamma_2$ is defined, in terms of the ``Gauss integral'', by:
 $$ \{\gamma_1,\gamma_2\} =\
\D{1\over 4\pi} \int_{t_1\in I_1}\int_{t_2\in I_2} \ \D {
<\g_1'(t_1)\wedge \g_2'(t_2), \g_{12}(t_1,t_2)> \over \Vert
\g_{12}(t_1,t_2)\Vert^3 } \ \ dt_1dt_2 , $$ where
$\g_{12}(t_1,t_2)=\g_2(t_2)-\g_1(t_1)$. \ev {\bf Remark.} In this
definition the normalization by $4\pi$ has to be seen as the
normalization by the volume of the unit sphere of $\d R^3$ (see
the Remark that follows the proof of Theorem $2.2$).
 \ev
Let us stress that this definition immediately shows  that
$\{\gamma_1,\gamma_2\}$ does not depend on the parametrizaton of
the curves nor on their rigid transformations. Intuitively
speaking the linking coefficient gives the algebraic (i.e. signed)
number of loops of one contour around the other. This
interpretation  is justified by the following result.

\ev
{\bf Theorem  2.2. --- }
Let $\g_1$ and $\g_2$ be two closed curves in $\d R^3$ and assume
that $I_1=[0,2\pi]$.
\ev

\item{(i)} {\sl The linking coefficient $\{\gamma_1,\gamma_2\}$ is an
integer, and is unchanged by deformations of
$\gamma _1$ and $\gamma_2$, involving no intersection of one curve
with the other.}

\item{(ii)} {\sl  Let $F:D^2\rightarrow \d R^3$  be a map of the disc
$D^2$ which agrees with  $\gamma_1:t\mapsto \g_1(t),\ 0\leq
t\leq 2\pi$, on the boundary $\partial D^2 \simeq S^1 \simeq
\D{[0,2\pi ]\over 0 \sim 2\pi} $, and is
transversal to the curve $\gamma_2\subset \d R^3$. Then the
``topological linking number" which is the intersection index
$F(D^2)\cdot \gamma_2$,  (i.e. the number of the intersection
points of $F(D^2)$ and $\gamma_2$,  counted with the signs,
reflecting the orientation), is equal to the linking coefficient
$\{\gamma_1,\gamma_2\}$.}
\item{}

\ev
{\bf Proof.}
 The closed curves
$\gamma_i(t),\ i=1,2$, give rise to a 2-dimensional,
closed, oriented parametric surface $\g_1\times \g_2$ in $\d R^6$:
$$ \gamma_1\times\gamma_2:(t_1,t_2)\mapsto(\g_1(t_1),\g_2(t_2)). $$
 Since the curves are non-intersecting the map
$\varphi:\gamma_1 \times\gamma_2\rightarrow S^2$, given by:
$$ \varphi(t_1,t_2)=\ {\g_1(t_1)-\g_2(t_2) \over
\Vert \g_1(t_1)-\g_2(t_2)\Vert}\, $$

\noindent is well defined. An easy geometric consideration shows
that the integrand in the Gauss integral  is just
the Jacobian of the map $\varphi$. Therefore the Gauss integral
above is equal to the degree of
the map $\varphi$. Hence the linking coefficient is indeed an
integer. Under deformations of the curves $\{\gamma_1,\gamma_2\}$
involving no intersection  one with the other, the map $\varphi$
undergoes a homotopy, so that its degree, and therefore also the
linking coefficient, are preserved. Let us stress that in the
process of these deformations each of the curves $\gamma_1, \
\gamma_2$ separately may cross itself in an arbitrary way.  Of course, the
topological linking number is also preserved by such
deformations.

We now prove (ii). If the curves are not linked (i.e. if by means
of a homotopy respecting non-intersection they can be brought to
opposite sides of a 2-dimensional plane in $\d R^3$) then
it can be verified directly that $\{\gamma_1,\gamma_2\}=deg\
\varphi=0$. In a general case we can ``push" $\gamma_1$ along
$\gamma_2$  in such a way that after this deformation it comes
close to $\gamma_2$ only in a neighborhood of exactly one point.
Then by applying another deformation (remind that
self-intersections of the curves are allowed) we reduce the
general case of the problem of calculating the linking coefficient
essentially to the following simple situation: the curve
$\gamma_2$ is a straight line, while $\gamma_1$ is a circle,
orthogonal to $\gamma_1$ and passed several times in the positive
or negative direction. Thus we suppose $\gamma_1$ and $\gamma_2$
to be given respectively by
$\g_1(t_1)=(cos\ t_1,\ sin\ t_1,0),\ 0\leq
t_1\leq 2\pi$ and
$\g_2(t_2)=(0,0,t_2),\
-\infty<t_2<\infty$. The linking coefficient for these two
curves is:
$$ \{\gamma_1,\gamma_2\}  =
{1\over 4\pi}\int^{\infty}_{-\infty} \int^{2\pi}_0
\D {dt_1  dt_2 \over (1+t^2_2)^{3/2}}
 ={1 \over 2} \int^{\infty}_{-\infty} \D{d t_2 \over (1+t_2^2)^{3/2}}$$
$$=
 {1\over 2} \int^{\infty}_{-\infty} {du \over ch^2( u)}=
{1\over 2}   [th (u)]^{+\infty}_{-\infty} = 1. $$
 Hence for these
two directed curves the statement (ii) of the Theorem holds. The
general result now follows via the deformation described above.\f

\ev
{\bf Remark.}
 One can give another proof of
Theorem $2.2$. As above, we notice that the Gauss integral is equal
to the degree of the mapping
$\varphi:\gamma_1\times\gamma_2\rightarrow S^2$. Now fix a point
$p$ in $S^2$ which is a regular value of $\varphi$ and consider
the projection $\pi$ along the corresponding line $\ell_p$ onto the
orthogonal plane $P_p$. The preimages $\varphi^{-1}(p)$ correspond
exactly to the crossing points of the plane curves $\pi(\gamma_1)$
and $\pi(\gamma_2)$ in $P_p$. The sign of the Jacobian of
$\varphi$ at each of the preimages $\varphi^{-1}(p)$ can be
computed via the directions of $\pi(\gamma_1)$ and $\pi(\gamma_2)$
at their corresponding crossing point, taking into account, which
curve is ``above" and which is ``below".

Now the degree of the mapping $\varphi$ is the sum of the signs of
the Jacobian of $\varphi$ over all the preimages
$\varphi^{-1}(p)$. On the other hand, the corresponding sum over
all the crossing points of $\pi(\gamma_1)$ and $\pi(\gamma_2)$ can
be easily interpreted as the topological linking number of
$\gamma_1$ and $\gamma_2$.
\ev
{\bf 2.2.2 Signed and absolute rotation}
\ev

As the curves $\gamma_1,\gamma_2$ are not necessarily closed, the
Gauss integral   can be still computed.
\ev
{\bf Definition. }
For two curves
$\gamma_1,\gamma_2$ in $\d R^3$, closed or non-closed, we
call the Gauss integral   along these curves the {\sl signed
rotation of the curves} $\gamma_1$ and $\gamma_2$ and denote it by
$R(\gamma_1,\gamma_2)$. We have:
$$ R(\gamma_1,\gamma_2) =\
{1\over 4\pi} \int_{t_1\in I_1}\int_{t_2\in I_2} \ {<\g'_1\wedge
\g'_2, \g_{12}>\over \Vert \g_{12}\Vert^3}\ \  dt_1dt_2 , $$
 where $\g_{12} =\g_2-\g_1$.

  For two curves $\gamma_1,\gamma_2$, not
necessarily closed, the {\sl absolute rotation} of these curves is
defined as:
 $$ R_{abs}(\gamma_1,\gamma_2)=\D{1\over 4\pi}
\int_{t_1\in I_1}\int_{t_2\in I_2} \ { \vert <\g'_1\wedge
\g'_2, \g_{12}>\vert \over \Vert \g_{12}\Vert^3}
\ \  dt_1dt_2 . $$
\ev
{\bf Remarks.}
The absolute rotation of two curves
bounds the absolute value of their signed rotation:
$$\vert R(\gamma_1,\gamma_2)\vert \leq
R_{abs}(\gamma_1,\gamma_2).$$ In particular, for
$\gamma_1,\gamma_2$ closed the absolute rotation bounds the
absolute value of the linking number $\{\gamma_1,\gamma_2\}$.

 For $\gamma_1,\gamma_2$ not closed,
$R(\gamma_1,\gamma_2)$ does not need to be an integer any more.
This rotation number also is not invariant under the deformations
of the curves $\gamma_1,\gamma_2$ without crossing one another,
even if we assume that the end-points of $\gamma_1$ and $\gamma_2$
are fixed. Indeed, if we take the curve $\gamma_1$ to be a ``long"
segment of the straight line, and the curve $\gamma_2$ to be the
unit circle around $\gamma_1$, the computation at the end of
the proof of Theorem $2.2$ above shows that the
signed rotation $R(\gamma_1,\gamma_2)$, is approximately one. On
the other hand, we can deform the circle $\gamma_2$, with one of
its points fixed, as follows: we pull it out from the segment
$\gamma_1$, and then contract it to the point. The rotation of the
deformed curves is zero, so it was not preserved in the process of
the deformation.

\ev

However, for one of the curves, say $\gamma_1$, closed, we have
the following result:
\ev
{\bf Proposition 2.3. --- }
{\sl Let the curve $\gamma_1$ be closed. Then the signed rotation
$R(\gamma_1,\gamma_2)$ is invariant under the deformations of the
curve $\gamma_2$ (without crossing $\g_1$) if the end-points of
$\gamma_2$ remain fixed.}

\ev {\bf Proof.} Consider a closed curve $\tilde \gamma_2$
obtained from $\gamma_2$ by passing it twice in the opposite
directions. We have $R(\gamma_1,\tilde \gamma_2)=0,$ since the
rotation of these two closed curves is invariant under
deformation, while $\tilde \gamma_2$ can be deformed into the
point without crossing $\gamma_1$. Now consider another
deformation of $\tilde \gamma_2$, where one copy of $\gamma_2$
remains fixed, while another copy, $\hat \gamma_2$, undergoes a
deformation without crossing $\gamma_1$ and with the end points
fixed. The signed rotation remain zero in this deformation, so
$R(\gamma_1,\tilde \gamma_2)= R(\gamma_1, \gamma_2)-
R(\gamma_1,\hat \gamma_2)=0.$ Hence $R(\gamma_1,\hat \gamma_2)$
remains the same in the deformation.\f \ev

Another property which can be obtained by a rather straightforward
computation, is the following:

\ev
{\bf Proposition 2.4. --- }
{\sl
For $\g_1=\r L$ a straight line,
the signed rotation $R(\gamma_1,\gamma_2)$
(resp. the absolute rotation $R_{abs}(\gamma_1,\gamma_2)$  )
given by the Gauss
integral coincides with the signed rotation $R(\g_2, \r L)$
(resp. the absolute rotation $R_{abs}(\g_2, \r L)$)
of the curve $\gamma_2$
around the straight line $\r L$, as defined by projection on
$\r L^\perp$ in Section $2.1$ above.

The absolute rotation $R_{abs}(\gamma_1,\gamma_2)$ given by the
expression (2.7) coincides with the absolute rotation of the curve
$\gamma_2$ around the straight line $\gamma_1$, as defined in
Section $2.1$. }

\ev
{\bf Proof. }
  Let us start with the case of the signed
rotation. We have, passing to the length parametrization,
$$
R(\gamma_1,\gamma_2)=\D {1\over 4\pi}
\int_{I_1}\int_{I_2} {< n_1\wedge n_2 ,\g_{12}>
\over \Vert \g_{12}\Vert^3} \ \ ds_1 ds_2,$$
 where $n_1,n_2$ are the
unit tangent vectors to the curves $\gamma_1,\gamma_2$, and $s_1,
s_2$ are the length parameters on the curves $\gamma_1,\gamma_2$,
respectively.

Let $\eta_2(s_2)$ denote the vector joining the point $s_2 \in
\gamma_2$ with the projection of $s_2$ onto the straight line
$\gamma_1$ (see Figure 1). We have $<n_1\wedge
n_2,\g_{12}>=<n_1\wedge n_2, \eta_2>$. Hence:
$$R(\gamma_1,\gamma_2)= \D{1\over 4 \pi} \int_{I_1}
\int_{I_2}{< n_1\wedge n_2 ,\eta_2>
\over \Vert \g_{12}\Vert^3} \ \ ds_1 ds_2. $$
 The vector
$n_1$, being the tangent vector to the straight line $\gamma_1$,
is constant. Therefore, the triple product under the above
integral depends only on $s_2$ and we have:
$$ R(\gamma_1,\gamma_2)=
\D {1 \over 4\pi} \int_{I_2} < n_1\wedge n_2(s_2) ,\eta_2(s_2)>
\int_{I_1} {ds_1 \over \Vert \g_{12}(s_1,s_2)\Vert^3} \ \ ds_2. $$
Now the integral: $$
\int_{I_1} {ds_1 \over \Vert \g_{12}(s_1,s_2)\Vert^3} $$
 over
the straight line $\gamma_1$ has been computed (up to a scaling by
the distance $\Vert \eta_2(s_2) \Vert$ to the line) at the end of
Section 2.2.1 above. It is equal to $\D {2\over {\Vert \eta_2(s_2)
\Vert}^2}$. So for the rotation $R(\gamma_1,\gamma_2)$ we finally
get:
$$ R(\gamma_1,\gamma_2)= {1\over {2\pi}}\int_{I_2} \D
{< n_1\wedge n_2(s_2) ,\eta_2(s_2)> \over {\Vert \eta_2(s_2)
\Vert}^2} \ \ ds_2. $$
 Now we can replace the vector $n_2$ in the triple
product by its projection $\hat n_2$ to the line orthogonal to the
lines $\gamma_1$ and $\eta_2$. Hence this triple product is equal
to $\D{{\Vert \hat n_2 \Vert}\over {\Vert \eta_2 \Vert}},$ with
the sign defined by the orientation, and we obtain:
$$ R(\gamma_1,\gamma_2)= {1\over
{2\pi}}\int_{\gamma_2}\pm \D {\Vert \hat n_2 \Vert\over \Vert
\eta_2 \Vert}\ \  ds_2 . \eqno(1)$$
 But $\eta_2$ is just the radius-vector of the
projection of the curve $\gamma_2$ onto the plane orthogonal to
the line $\gamma_1$, and $\hat n_2$ is the orthogonal component of
the velocity vector of this projection. Therefore, according to
Lemma $2.1$ and up to a sign $\epsilon$,
 the integral $(1)$ is the signed rotation $R(\g_2,\r L)$ of
$\gamma_2$ around the straight line $\r L=\gamma_1$, as defined in
Section $2.1$ (the sign $\epsilon$ is $+$ when in the computation
of $R(\g_2,\r L)$ we have oriented the plane $\r L^\perp$ in such
a way  that if $(\vec u,\vec v)$ is an oriented orthonormal basis
of $\r L^\perp$, $\vec u\wedge \vec v=-n_1$ (see Figure $1$)).
This completes the proof of the Proposition for the signed
rotation.

\vskip7.5cm
\vglue0cm
\vskip0cm
\hskip0,5cm\includegraphics{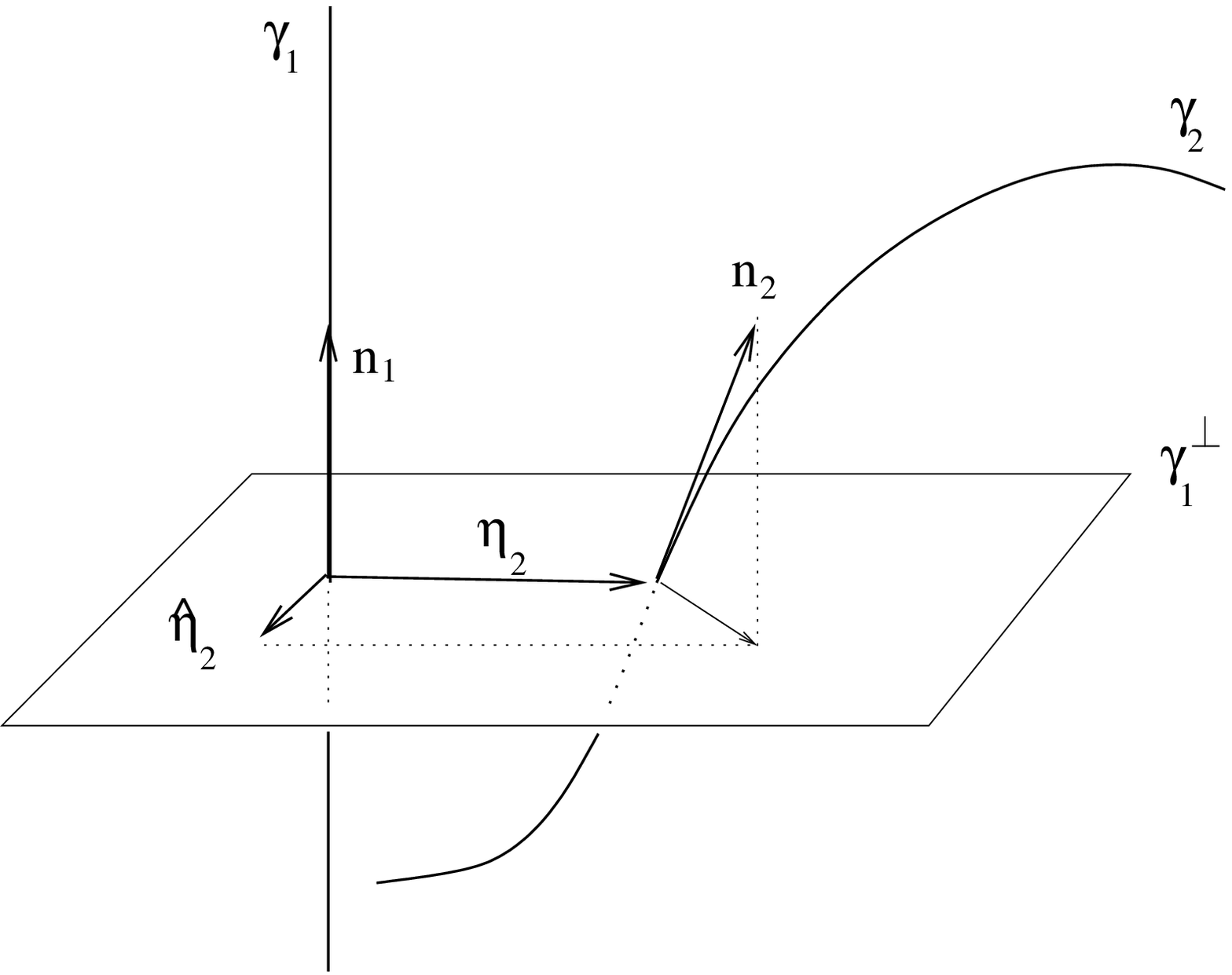}
\vglue0cm
\Ev
\centerline{\bf Figure 1}
\ev

The proof for the absolute rotation is exactly the same: we just
take the absolute value of the triple product in each step. Since
the proof above consisted of a chain of point-wise equalities of
the integrands, it remains valid also for the absolute values.
\f

\vfill\eject

\centerline{\bf 3. Rotation of trajectories of Lipschitz vector fields}
\Ev For a Lipschitz vector field in $\d R^n$, the very simple and
basic fact is that the angular velocity of its trajectories with
respect to any {\sl stationary} point is bounded by the Lipschitz
constant $K$. As an immediate consequence, the rotation of a
trajectory around a stationary point is bounded (up to some
universal constant $C$) by the Lipschitz constant $K$ and the time
interval $T$. The bound has this form: $C\cdot K\cdot T$.

Below we remind the proof of this fact. In fact, we show
that the same is true for the rotation speed of the trajectories
around any linear subspace ${\cal L} \subset \d R^n$, which
is an invariant submanifold of our field.

\ev

The main result of this section (Theorem $3.4$) is that the
``rotation" of a trajectory around any point ({\sl stationary or
not}) is bounded in terms of $K$ and $T$. As the consequence we
prove (Theorem $3.8$) that the mutual rotation of any two
trajectories of a Lipschitz vector field is essentially bounded by
the Lipschitz constant and the length of the time interval. More
accurately, the absolute rotation of any two trajectories of a
Lipschitz vector field on the time intervals $T_1$ and $T_2$,
respectively, is bounded from above by a linear combination of the
expressions $K\cdot \min(T_1, T_2)$ and $C \cdot K^2 \cdot T_1
\cdot T_2$. Easy examples, given in the end,  show that this bound
is sharp. \ev At this point, let us remind that the rotation of a
curve $\omega:I\to \d R^n $ around a point, an affine space or
another curve, as defined above,  only depends on the trajectory
$\omega(I)$, provided we only admit injective (parametrization of)
curves. Consequently our bounds not really depend  on the field,
but rather on the geometry of the trajectories of the fields. This
is obvious when we look at the type of bounds we obtain: (up to
some universal constant) they are combinations of $K\cdot T $ and
$K^2\cdot T_1\cdot  T_2$, expressions that are unchanged with
respect to any transformation of the field that preserve the
trajectories. \ev
 For Lipschiz vector fields, we can summarize the situation by the
following slogan: {\sl``two trajectories have finite mutual rotation in
finite time".}
\ev

We cannot expect bounds of this sort to be true for a rotation
around a non-invariant subspace of the field. Indeed, the
easy-to-construct Example
$3.3$ below shows that a trajectory of a $C^{\infty}$-vector field
in $\d R^3$ can make in {\sl finite} time an {\sl infinite} number of
turns around a straight line.

\ev

Let us remind that it was shown (in some special cases) in
[Gr-Yo] that for a trajectory of a
{\sl polynomial} vector field (trajectory which is not in general
in some o-minimal structure), its rotation rate around any algebraic
submanifold is bounded in terms of the degree of the submanifold
and the degree and size of the vector field. As a consequence, we
obtain also a linear in time  bound on the number of intersections
of the trajectory with any algebraic hypersurface.

On the other hand, our bounds on the rotation rate of the
trajectories of a polynomial vector field imply upper bounds on
the multiplicities of the local intersections of such trajectories
with algebraic submanifolds in terms of the degree only.

\ev
\ev

As the Example $3.3$ shows, nothing of this sort can be expected
even for $C^{\infty}$ (and of course, for Lipschitz) vector
fields. Still, the results of this section show that Lipschitz
vector fields exhibit rather strong non-oscillation patterns. As
far as the rotation of the trajectories of such fields around
non-invariant submanifolds is concerned, our current understanding
is far from being sufficient. In particular, there is a serious gap
between the result of Theorem $3.4$ below that a ``global"
rotation rate of a trajectory of a Lipschitz vector field around a
{\it non-stationary point} is still bounded, and the Example~$3.3$,
demonstrating an infinite rotation around a non-invariant straight
line.
\Ev
{\bf 3.1. Rotation of a trajectory around a stationary point}
\ev
Let  $v$ be a vector field defined in a certain domain $U$ in
$\d R^n$. We shall always assume $v$ to satisfy a Lipschitz
condition with the constant $K$:
$$ \Vert v(x)-v(y)\Vert \leq
K\Vert x-y\Vert, \ \hbox{ for any two points }\  x,y \in U$$

\ev

\noindent  Let $x_0 \in U$ be a stationary or a singular point of
$v$, ie $v(x_0)=0$, so that the constant curve $c(t)=x_0, t\in \d
R,$ is the integral curve of $v$ passing through $x_0$. Then for
any $x \in U, \ x \ne x_0$, the angular velocity of the trajectory
of $v$, passing through $x$, with respect to $x_0$, is equal to
${\Vert \hat v(x)\Vert} / {\Vert x-x_0\Vert}$, where $\hat v(x)$
is the projection of the vector $v(x)$ to the hyperplane
orthogonal to $x-x_0$. Hence, this angular velocity does not
exceed $K$:

$$ {{\Vert \hat v(x)\Vert} \over
{\Vert x-x_0\Vert}} \leq {\Vert v(x)\Vert\over \Vert
x-x_0\Vert}={\Vert v(x)-v(x_0)\Vert \over\Vert x-x_0\Vert}\leq K
$$
 where $K$ is the Lipschitz constant of $v$. By Lemma $2.1$ we
obtain :

\ev
{\bf Proposition 3.1. --- }
{\sl For any trajectory $\omega(t)$ of the field $v$, its
rotation around the stationary point $x_0$ of the field between
the time moments $t_1$ and $t_2$, i.e., the length of the
spherical curve $s(t)=\D {\omega(t)-x_0\over\Vert \omega(t)-x_0\Vert}$ between
$t_1$ and $t_2$, does not exceed $K\cdot (t_2-t_1)$.}
\ev
{\bf Remark. }
Of course on a non-finite time interval, a trajectory may have
infinite local rotation  around a stationary point, as shown by
the example given below (we avoid the obvious example of a cyclic
trajectory, since we aim to work with injective trajectories).

Let us consider the following
algebraic field in $\d R^2$ with singular point $O=(0,0)$:
$$v(x,y)=\Big( (x^2+y^2-1)x-y, (x^2+y^2-1)y+x \Big),$$
and introduce the following notations:
 $r^2= x^2+y^2$, $\vec\tau=(x,y)/r$ and $\vec n =(-y,x)/r$.
 We have :
 $$v(x,y)=r(r^2-1)\cdot \vec \tau+r \cdot \vec n. $$

A trajectory passing at a point $p$ with $r(p)=1$ has to be
 the unit circle. Now consider $\omega=(\alpha, \beta)$ an integral
curve of $v$ passing through a point $q$ with $r(q)<1$. This
trajectory cannot go outside the open unit disc and as we have:
$\D{d[r^2\circ \omega]\over dt}=2(\alpha\alpha'+\beta\beta')$, we
obtain: $\D{d[r^2\circ \omega]\over dt}=
2(r^2-1)(\alpha^2+\beta^2)<0$. This proves that the trajectory
$\omega$ has the singular point $O$ as limit (while the unit
circle has to be a limit cycle of this trajectory). Furthermore we
know that from the point $q$, the limit $O$  is approached on an
infinite time interval $I$.

On the other hand, the velocity $\ti v$ of the spherical blowing-up
$\sigma$
of $\omega$ is: $\D{1\over r}\cdot r\cdot \vec n=\vec n$.
We conclude that :
$$R_{abs}(\omega,O)=\D \int_I \Vert \vec n(w(t))\Vert\ dt =
\int_I 1 \ dt=+\infty.$$

\vskip7.3cm
\vglue0cm
\vskip0cm
\hskip1.6cm\includegraphics{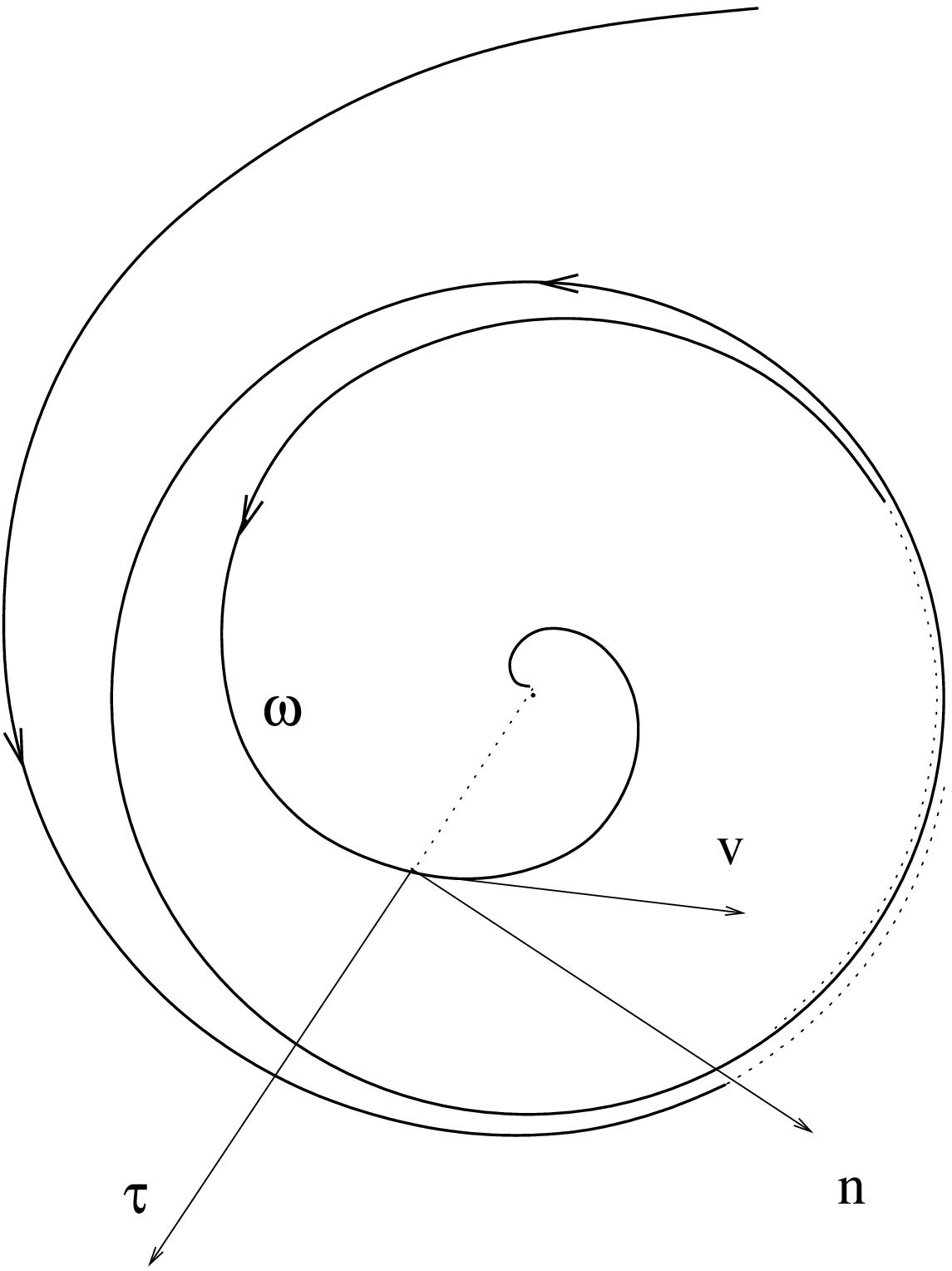}
\vglue0cm
\ev
\centerline{\bf Figure 2}
\Ev

To finish with this remark, let us recall that such an example of
spiraling trajectory does
not exist for $v$ a gradient vector field of an
analytic map, as proved in [Ku-Mo-Pa], [Ku-Mo]; for such a field
the length of $\sigma$ is finite and the limit
of the secants passing through the singular point does exist
 (Thom's Gradient Conjecture).

\ev
 Exactly in the same way as in Proposition $3.1$
we prove the following more general result
about rotation of a curve around an invariant affine space $\r L$.
The proof actually shows that in Proposition~$3.1$ the assumption
$v(x_0)=0$ has preferably to be considered  as ``{\sl $v(x_0)$ is
tangent
to the submanifold $\{x_0\}$ of $\d R^n$}", or ``{\sl
$v$ is a stratified
field with respect to the stratification $(\{x_0\}, U\setminus
\{x_0\})$ of $U$}".

 \ev
{\bf Proposition 3.2. --- }
 {\sl For any affine subspace ${\cal L} \subset \ \d R^n$
which is invariant for the vector field $v$ (i.e. for any $x\in
{\cal L}$, $v(x)$ is tangent to ${\cal L}$ or, in other words, $v$
is stratified with respect to $(\r L, U\setminus \r L)$), the
rotation speed of $v$ in the orthogonal to ${\cal L}$ direction is
bounded by $K$. In particular, the absolute rotation of any
trajectory $\omega$ of the field $v$ around $\cal L$ in time $T$
does not exceed $K\cdot T$.}

 \ev
{\bf Proof.}
According to the definition of the absolute rotation
around a linear subspace (Section $2.1$), we consider the
orthogonal to $\cal L$ component $\tilde v(x)$ of the vector field
$v(x)$. Let $x_0$ be the projection of $x$ onto $\cal L$. We
denote by $\hat v(x)$ the ``rotation" component of $\tilde v(x)$,
orthogonal to $x-x_0$. Then the rotation speed of $v$ in the
orthogonal to $\cal L$ direction is equal to ${\Vert {\hat
v}(x)\Vert}/{\Vert x-x_0\Vert}$. Hence, this rotation speed is
bounded as follows:
$$  {\Vert {\hat v}(x)\Vert \over \Vert
x-x_0\Vert} \leq {\Vert \tilde v(x)\Vert\over \Vert x-x_0\Vert}=
{\Vert \tilde v(x)- \tilde v(x_0)\Vert \over\Vert x-x_0\Vert}\leq
{\Vert v(x)- v(x_0)\Vert \over\Vert x-x_0\Vert}\leq
K. $$
Here we use the fact that $\cal L$ is invariant for $v$ and
hence $\tilde v(x_0)=0$. By Definition  we obtain that the
absolute rotation of any trajectory $w(t)$ of the field $v$ around
$\cal L$ in time $T$ does not exceed $K\cdot T$. This completes the
proof of Proposition $3.2$.\f

\ev
If $x_0$ is not a stationary point of $v$, the
angular velocity of $v(x)$ with respect to $x_0$ tends to infinity
as $x$ approaches $x_0$ in any  direction transverse to $v(x_0)$.
Indeed, as $x$ approaches $x_0$, $v(x)$ tends to $v(x_0)\ne 0$. By
Lemma $2.1$, the angular velocity of $v(x)$ with respect to $x_0$ is
equal to ${\Vert {\hat v}(x)\Vert}/{\Vert x-x_0\Vert}$, where
$\hat v (x)$ is the component of $v(x)$ orthogonal to $x-x_0$. So
if $x$ approaches $x_0$ in a direction transverse to $v(x_0)$,
$\hat v(x)$ tends to $\hat v_0\not=0$, while $x-x_0$ tends to zero. See
Figure~$3$.

\vglue0cm
\vskip4.5cm
\vglue0cm
\vskip0cm
\hskip1,5cm\includegraphics{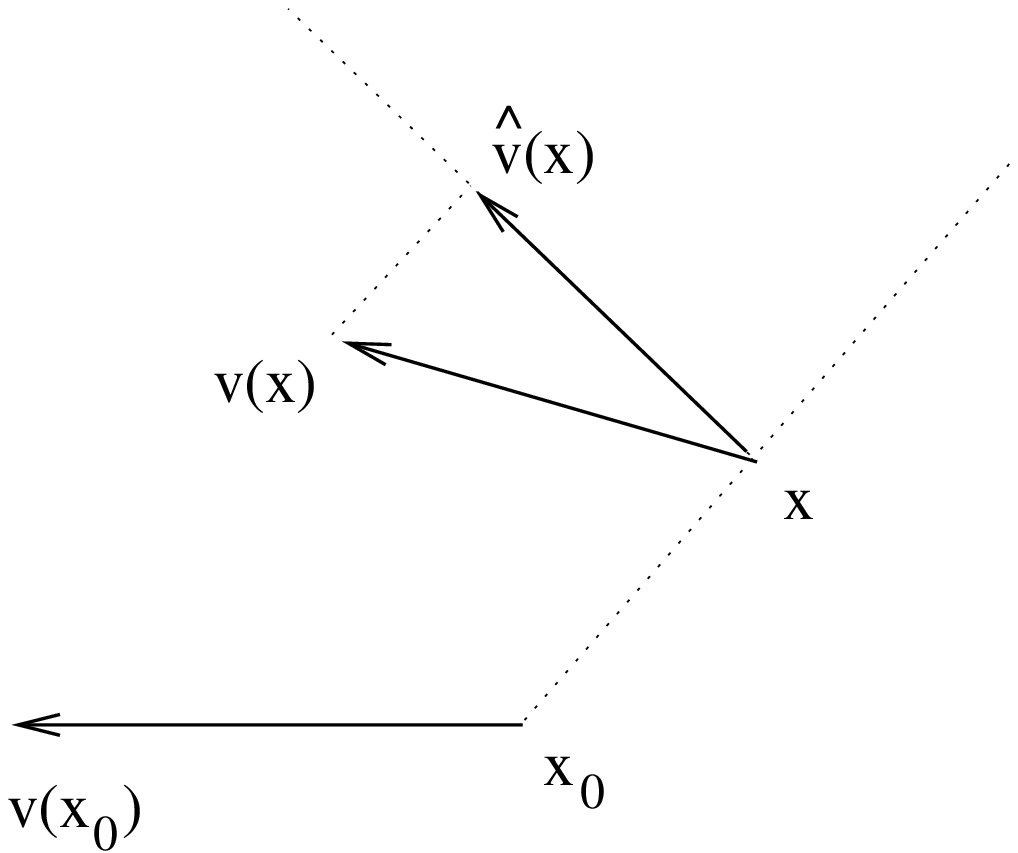}
\vglue0cm
\Ev
\centerline{\bf Figure 3}
\Ev

However, one can show (see Theorem $3.4$ below) that for any
trajectory $\omega(t)$ of $v$, and for $t_2-t_1$ big enough, the
length of the spherical curve $s(t)={\omega(t)-x_0/\Vert
\omega(t)-x_0\Vert}$ (i.e. the absolute rotation of $\omega(t)$
around $x_0$) is still bounded by $C \cdot K\cdot(t_2-t_1)$. The
intuitive explanation is as follows: as the trajectory $\omega(t)$
passes very close to $x_0$, its rotation around $x_0$ in the time
interval $[t_1,t_2]$ comes to approximately $1/2$ (see Figure
$4$).

\vskip4.5cm
\vglue0cm
\vskip0cm
\hskip0,8cm\includegraphics{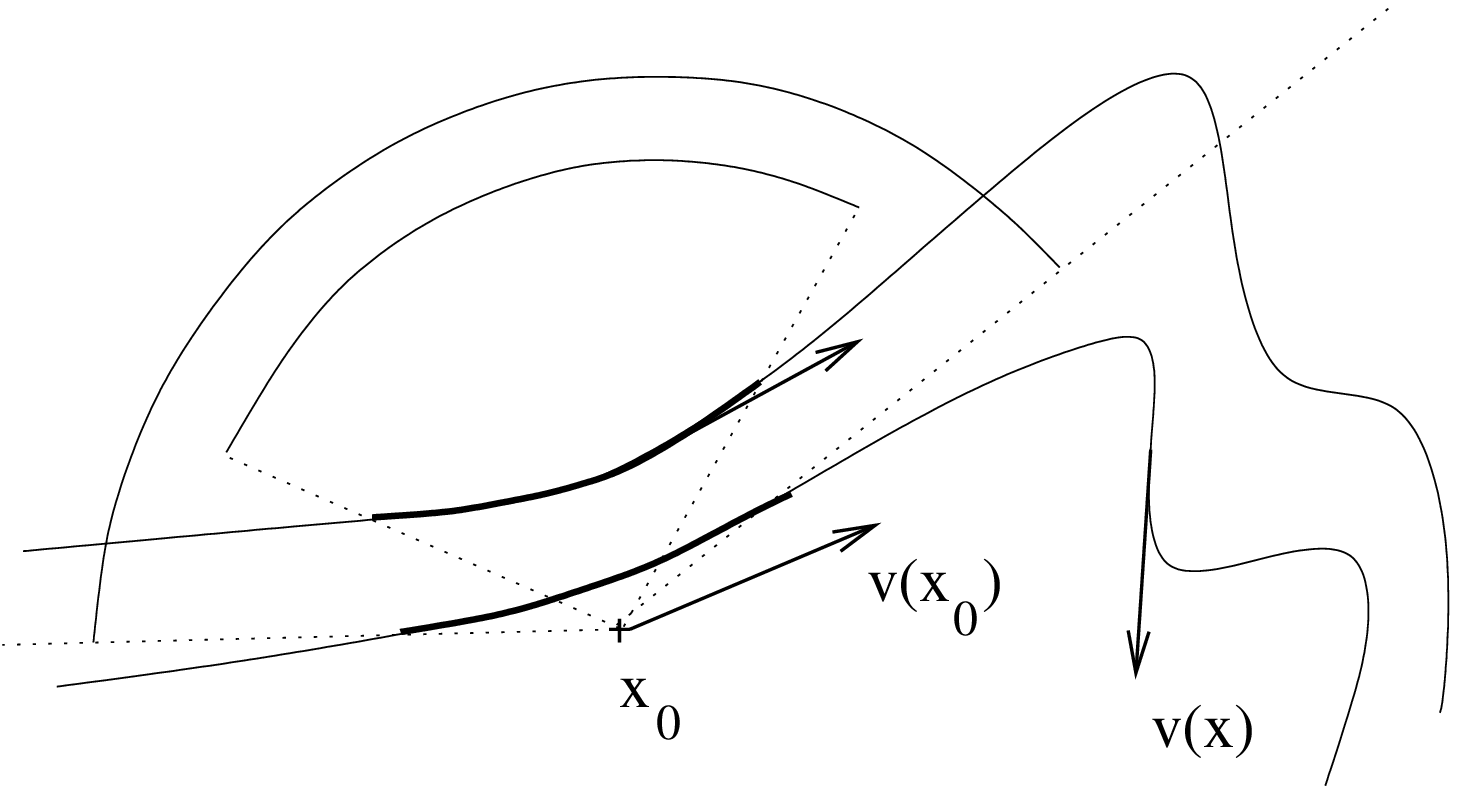}
\vglue0cm
\ev
\centerline{\bf Figure 4}
\Ev
 However, continuing to move in the same direction, the
trajectory $\omega(t)$ cannot gain more rotation. So for the total
rotation of the trajectory $\omega(t)$ around $x_0$ to grow, its
velocity vector $v(\omega(t))$ must change its direction. On the
other side, if we know a priori that {\it the directions of $v(x)$
and $v(x_0)$ are strongly different with respect to $K\cdot \Vert
x-x_0\Vert$}, then the rotation speed satisfies the same Lipschitz
upper bound as above. Consequently, the time interval $[t_1,t_2]$
has to be large in order to gain a large total rotation.

\ev 

Although the rotation in finite time of a trajectory of a
Lipschitz field around a point is bounded, the rotation of a
trajectory   of a Lipschitz vector field $v$ around a straight
line, {\sl which is not invariant under $v$}, can be unbounded
during a finite interval of time,
 and this phenomena may occur even
for $v$ a $C^{\infty}$-vector field.
Consider for instance the following field:

\ev

\noindent
{\bf Example 3.3}.
 Let $\Phi: \d R^3\to
\d R^3$ be a diffeomorphism, defined by:
$$\Phi(x_1,x_2,x_3)=(x_1,x_2,x_3),\ \ x_1\leq 0,$$
$$\Phi(x_1,x_2,x_3)=(x_1,x_2+\omega_1(x_1),x_3+\omega_2(x_1)),\
\ x_1\geq 0,$$
where $\omega_1(x_1) = e^{-1/x^2_1}\cos({1/ x_1})$, and
$\omega_2(x_1) = e^{-1/x^2_1}\sin({1/ x_1})$.

\ev

\noindent One can easily check that $\Phi$ is a
$C^{\infty}$-diffeomorphism of a neighborhood of $0\in \d R^3$.
Now the image of the positive $x_1$-semiaxis under $\Phi$ is
a line $w$, which makes an infinite number of turns around $Ox_1$
in any neighborhood of the origin.

\ev

Consider the vector field $v$ in $\d R^3$, which is an image under
$\Phi$ of the constant vector field $e_1=(1,0,0)$. Clearly, $w$ is a
trajectory of the $C^{\infty}$-vector field $v$, and it makes an
infinite number of turns around $Ox_1$ in finite time. In coordinates,
$$v(x_1,x_2,x_3)=D\Phi_{(x_1,x_2,x_3)}(e_1)=\D{\partial \Phi \over
\partial x_1}(x_1,x_2,x_3)=(1,{\omega}'_1(x_1),{\omega }'_2(x_1)).$$
Notice that in this example, the orthogonal components of $v$ on
the line $Ox_1$ itself has an infinite number of sign changes,
accumulating to the origin.

\ev

\noindent
{\bf Remark.}
Notice that the rotation of any two
trajectories of the vector field of the Example $3.3$ one around
another is zero. Indeed, the field $v$ does not depend on the
coordinates $x_2, \ x_3$. So the vector joining the intersection
points of the two trajectories with the planes $x_1=c$ remains
constant.
In particular, this shows that we cannot expect any ``transitivity"
in the rotation of three curves: take two trajectories of $v$
``far away" from the line $Ox_1$, while the third trajectory is
$w$ as in Example $3.3$.
\Ev
\ev
{\bf 3.2. Rotation of a trajectory around a non-stationary point}
\ev

As it was mentioned above, one cannot
bound uniformly the {\it momentary angular velocity} of trajectories of a vector field $v$
with respect to a {\sl non-stationary} point $x_0$.
We shall show in this section that nevertheless the ``long-time"
rotation rate of trajectories of a Lipschitz field $v$ {\it with
respect to any point} $x_0$, stationary or non-stationary, is
uniformly bounded.

\ev

Our main result is the following Theorem: \ev
{\bf Theorem 3.4. --- }
{\sl Let $v$ be a Lipschitz vector field on an open set $U \subset
{\d R}^n$ with Lipschitz constant $K$. Let $\omega(t)$ be a
trajectory of $v$. Then for any $x_0 \in U$ the absolute rotation
of $\omega$ around $x_0$ between any two time moments $t_1$ and $t_2$
satisfies:
 $$R_{abs}(\omega,x_0)\le 4+
 K\cdot (t_2-t_1).$$ }

\ev
{\bf Proof.}
To prove Theorem $3.4$ we need a general geometric lemma, which
expresses accurately (in one of many possible ways) the
intuitively clear fact that a long curve inside a fixed sphere
must oscillate. Let $\omega(t)$ be a $\r C^1$-curve (or piecewise
$\r C^1$) in the unit sphere $S^{n-1} \subset \d R^n$. For any two
time moments $t_1$ and $t_2$ denote by $s(t_1,t_2)$ the length of
the curve $\omega$ between $t_1$ and $t_2$.

For each ``equator circle" $\lambda$ in $S^{n-1}$ consider the
``longitude" projection $\pi: S^{n-1} \rightarrow \lambda$ of the
sphere (without the poles) on the circle $\lambda$. A circle
$\lambda$ is given as the intersection of $S^{n-1}$ with a
$2$-dimensional vector plane $\Lambda$ of $\d R^n$. The fiber of
$\pi $ over a point $t\in \lambda$ lies in a $(n-2)$-dimensional
sphere $S^{n-2}_{\lambda,t}\subset S^{n-1}$ which is normal to
$\lambda$ and characterized as  the intersection $H\cap S^{n-1}$,
where $H$ is the hyperplane of $\d R^n$ passing through $t$ and
containing $\Lambda^\perp$. The projection of a point $x\in
S^{n-1}$ is defined as soon as $x$ is not in the poles, that is to
say $x\not \in \Lambda^\perp\cap S^{n-1}\simeq S^{n-3}$(see Figure
$5$). Since $\omega$ is (piecewise) $\r C^1$, by Sard's theorem,
it does not intersect a generic sub-sphere $S^{n-3}\subset
S^{n-1}$. Hence, the projection $\pi:\omega \to \lambda$ is
defined and regular for generic circles $\lambda$.

For the computations below it is convenient to fix a real
parameter $\theta > 4$.\ev
{\bf Lemma 3.5. --- }
{\sl Let $R>0$ and $S^{n-1}_R$ be the sphere of $\d R^n$ of radius
$R$. If the length $s(t_1,t_2)$ of a curve $\omega\subset
S^{n-1}_R$ between $t_1$ and $t_2$ satisfies $s(t_1,t_2)>2\pi
R\cdot \theta$ then there exist an ``equator circle" $\lambda$ in
$S^{n-1}$ and time moments $\tau_1$ and $\tau_2$, with $ t_1 <
\tau_1 < \tau_2 < t_2$, such that the following is true:}

\ev

\item{(i)} {\sl The ``longitude" projections $\pi(\omega(\tau_1))$ and
$\pi(\omega (\tau_2))$ of the points $\omega (\tau_1)$ and
$\omega(\tau_2)$ on the circle $\lambda$ coincide or are antipodal
points on $\lambda$. }

\item{(ii)} {\sl Denoting  by $\ell$ the tangent line to
$\lambda$ at this common projection point or antipodal points,
the norm of the projections of the velocity vectors
$v(\tau_i)=
 \D{d\omega \over dt} (\tau_i), \ i=1,2$ on the line $\ell$ are
$\ge \D{\theta - 4 \over 4\theta} \cdot \D {s(t_1,t_2)\over
t_2-t_1}$, and the directions of these projections are opposite in
$\ell$ (in case $\pi(\omega(\tau_1))$ and $ \pi(\omega(\tau_2))$
are not equal but antipodal, $v(\tau_2)$ is identified with its
projection on the tangent line at $\pi(\omega(\tau_1))$). } \ev
We prove
Lemma $3.5$ below in this section. Let us complete now the proof of
Theorem~$3.4$.

\ev
We apply Lemma $3.5$ to the spherical curve
$\sigma(t)=(\omega(t)-x_0)/\Vert \omega(t)-x_0\Vert$, which is
contained in
the unit sphere $S^{n-1}$ in $\d R^n$. By this lemma either
the length $s(t_1,t_2)$ of $\sigma(t)$ between the moments $t_1$
and $t_2$ is smaller than $2\pi\cdot \theta$, or there is an equator
circle $\lambda \subset S^{n-1}$ and time moments $\tau_1, \tau_2$
such that the longitude projections $\pi(\sigma_1)$ and
$\pi(\sigma_2)$ of the points $\sigma_1=\sigma(\tau_1),
\sigma_2=\sigma(\tau_2)$ coincide or are antipodal points on $\lambda$.
 For instance, let us assume that $\pi(\sigma_1)=\pi(\sigma_2)=\mu$,
the proof being the same in the case $\pi(\sigma_1)$ and
$\pi(\sigma_2)$ are antipodal in $\lambda$, because the tangent lines
to $\lambda$ at antipodal points are the same.
  By Lemma $3.5$, the orthogonal projections of the velocity vectors
of $\sigma$ at $\tau_i, \ i=1,2$ on $\ell$, the tangent line  to
$\lambda$ at $\mu$, are  in
norm  $\ge \D {\theta-4 \over 4\theta}\cdot{s(t_1,t_2)
\over t_2-t_1}$ and
the directions of these two projections are opposite.

\ev

Denoting the corresponding points $\omega(\tau_1)$ and
$\omega(\tau_2)$ of the original trajectory $\omega(t)$ by $w_1$
and $w_2$, we have: $\D\sigma_1={w_1-x_0\over \Vert w_1 -x_0
\Vert} $ and $ \D \sigma_2={w_2-x_0 \over\Vert w_2-x_0 \Vert}$.
   Let $v(t)$ denote, as above,
the velocity vector of the trajectory $\omega(t)$, while $\tilde v(t)$
denotes the velocity vector of the spherical curve $\sigma$. The
vector $\tilde v(t)$ is the orthogonal to $\omega(t)-x_0$ component of
$\D{v(t)\over \Vert \omega(t)-x_0 \Vert}$.

\ev

A simple property of the longitude projection is the following
(see Figure $5$): if for a point $s \in S^{n-1}$ we denote by
$T(s)$ the tangent hyperplane to $S^{n-1}$ at $s$ (which is also
the orthogonal hyperplane to the radius -vector of $s$), and by
$\ell(\pi(s))$ the tangent line to the circle $\lambda$ at the
longitude projection $\pi(s)$ of $s$ to $\lambda$, then $\ell(\pi(s))
\subset T(s)$. In particular, denoting by $T_1,T_2$ the tangent
hyperplanes to $S^{n-1}$ at $\sigma_1$ and $\sigma_2$,
respectively, we have $\ell \subset T_1, \ \ell \subset T_2.$

\vskip6,5cm
\vglue0cm
\vskip0cm
\hskip0cm\includegraphics{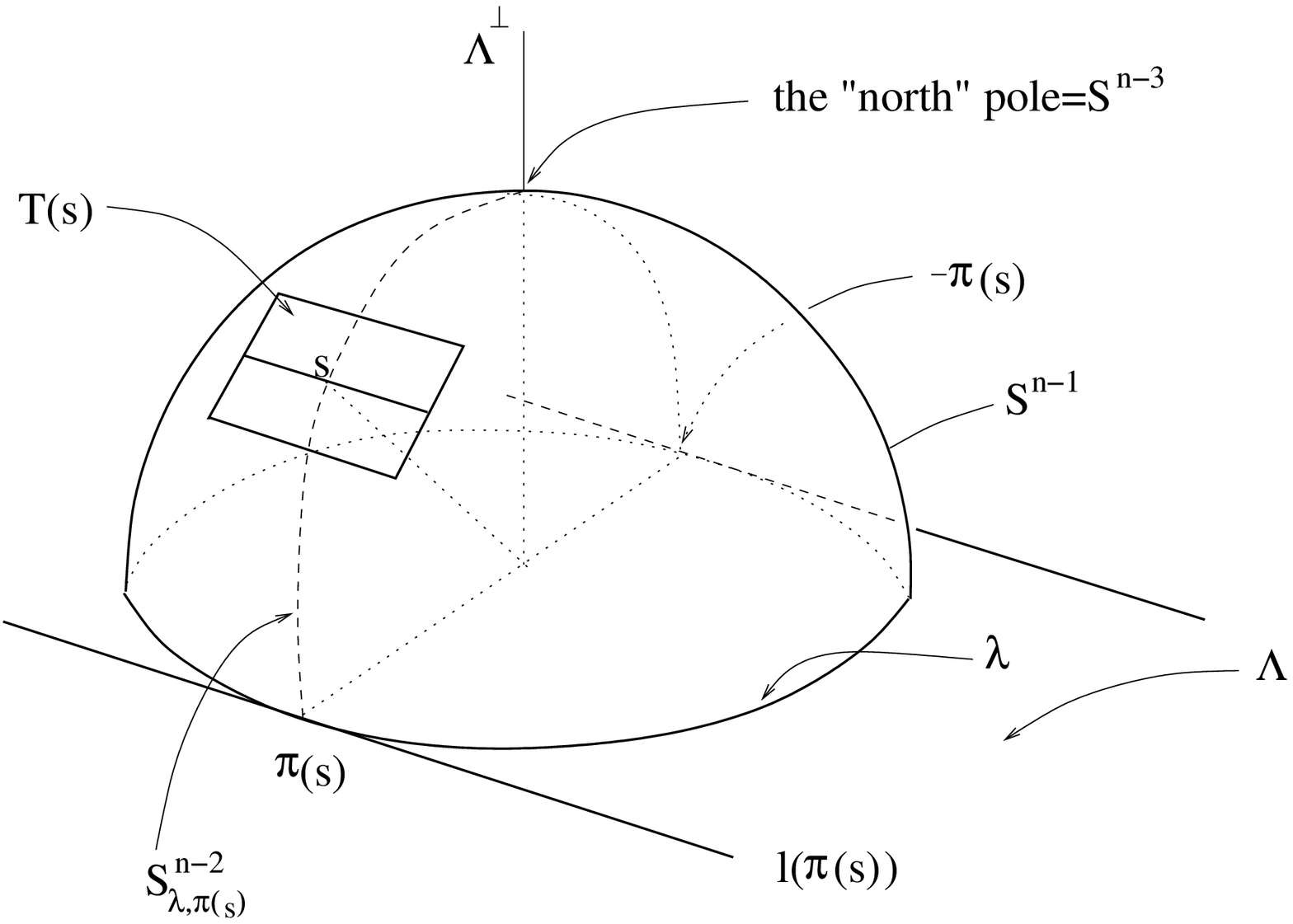}
\vglue0cm
\ev
\centerline{\bf Figure 5}
\ev
The immediate consequence of this inclusion is that the orthogonal
projections on $\ell$ of the velocity vectors $\tilde v_1 $ and
$\tilde v_2$ of the spherical
curve $\sigma$ are the same as the projections on $\ell$ of the
normalized velocity vectors $\D{v_1 \over \Vert w_1-x_0 \Vert}$
and $\D{v_2 \over \Vert w_2-x_0 \Vert}$ of the original
trajectory~$\omega$. Consequently, these two projections are
 $\ge\D {\theta-4\over 4\theta} \cdot {s(t_1,t_2)\over t_2-t_1}$ and their
directions are opposite.

\ev

Let us now consider two auxiliary vectors: $V_1=\D {v_1-v(x_0)
 \over \Vert w_1-x_0 \Vert}$ and $V_2=\D{{v_2-v(x_0) \over \Vert
w_2-x_0 \Vert}}$.
 Since the projections of $\D {v_1 \over \Vert
w_1-x_0 \Vert}$ and $\D{v_2 \over \Vert w_2-x_0 \Vert}$ on $\ell$
had opposite directions, independently of the vector $v(x_0)$ at
least one of the vectors $V_1,V_2$ has its projection on $\ell$ of
 size  $\ge \D {\theta - 4 \over 4\theta}\cdot
{s(t_1,t_2)\over t_2-t_1}$. But since the vector field $v$ is
Lipschitz, the norms of both the vectors $V_1$ and $V_2$ are
bounded by $K$. Hence $\D{\theta - 4\over 4\theta} \cdot
{s(t_2,t_1)\over t_2-t_1}\leq K$, or
$\D s(t_2,t_1)\leq {4\theta\over \theta-4}\cdot K\cdot (t_2-t_1)$.
In any case we have proved :
$$ R_{abs}(\omega, x_0)=\D
s(t_2,t_1) \leq \max(\theta , {4\theta\over (\theta-4)}\cdot
K\cdot (t_2-t_1) ). $$ But this bound is true for any $\theta>4$,
and for $\theta=4+ K\cdot (t_2-t_1)$, we have
both expressions under the maximum sign equal one another.
Therefore, for this specific choice of $\theta>4$ we get:
$$ \max(\theta ,
{4\theta\over (\theta-4)}\cdot K\cdot (t_2-t_1) )= 4+
 K\cdot (t_2-t_1)$$ Up to the proof of Lemma $3.5$,  this
completes the proof of Theorem $3.4$.\f

\ev
{\bf Remark.}
We could complete the proof of Theorem $3.4$
without returning to the point $x_0$, just by comparing the
vectors $v_1$ and $v_2$ and using the Lipschitz property of the
vector field $v$. However, the comparison with the vector $v(x_0)$
used in the proof above, can be applied also in more general
situations, like estimating rotation of a Lipschitz vector field
around a non-invariant subspace.

\ev

\ev
{\bf Proof of Lemma 3.5.}
Let $\omega(t)$ be a $\r C^1$-curve (possibly piecewise $\r C^1$)
in $S^{n-1}_R$ and assume that $\omega$ is defined on the time
interval $[t_1,t_2]$. We denote, as above, by $s(t_1,t_2)$ the
length of $\omega$ (between $t_1$ and $t_2$). For each one
dimensional circle $\lambda \simeq  S^1_R$ in $ S^{n-1}_R $,
 consider the longitude projection
$\pi_\lambda: S^{n-1}_R\rightarrow \lambda$. For $t\in \lambda$ we
denote by $S^{n-2}_{\lambda,t}$ the $(n-2)$-sphere of $S^{n-1}_R$
which contains the fiber $\pi_\lambda^{-1}(t)$. We denote by $\r
O(n)$ the orthogonal group of $\d R^n$ and by $\theta_n$ its Haar
measure. For $S \simeq S^1_R$ a fixed circle in $S^{n-1}_R$, $g\in
\r O(n)$ and $t\in S$, we denote by $S^{n-2}_{g\cdot S,g\cdot t}$
the $(n-2)$-sphere of $S^{n-1}_{R}$, which contains the fiber over
$g\cdot t$ of the longitude projection onto $g\cdot S$. By the
spherical Cauchy-Crofton formula (see [Fe 2] $3.2.48$), we have:

$$ \D{1\over \pi R}\cdot  s(t_1,t_2) = \int _{g\in {\cal O}(n)}
\#(\omega\cap S^{n-2}_{g\cdot S, g\cdot t})\ d\theta_n(g)$$
Now by Fubini's theorem:
$$  2 \cdot s(t_1,t_2) = \int _{t\in S} \
\int _{g\in {\cal O}(n)}
\#(\omega\cap S^{n-2}_{g\cdot S, g\cdot t})\ d\theta_n(g)\ dt$$
$$=
\int _{g\in {\cal O}(n)} \ \int _{t\in S} \
\#(\omega\cap S^{n-2}_{g\cdot S, g\cdot t})\ dt\  d\theta_n(g)$$
We conclude that for some $g\in \r O(n)$,
$$
\int _{t\in S} \
\#(\omega\cap S^{n-2}_{g\cdot S, g\cdot t})\ dt\ge 2\cdot s(t_1, t_2) $$
  But this integral is exactly twice
the length of the projection of $\omega$ onto the circle
$\lambda=g\cdot S$.

\ev

So we have shown that for any curve $\omega  \subset \d R^n$
there is a circle $\lambda\in S^{n-1}_R$ such that the
projection of $\omega$ onto $\lambda$ has its length
$\ge s(t_1,t_2) $.
But by our assumption, we also have:
$$ s(t_1,t_2)> 2\pi R\cdot \theta$$
Consequently, to  prove Lemma $3.5$ it is enough to prove
the following Lemma $3.6$, which is Lemma $3.5$  in
the case of a curve $\omega$ contained in a circle $S^1_R$:
\vfill\eject
\Ev
{ \bf Lemma 3.6. --- }
{\sl Let $\theta$ be any real number $>4$ and $\omega$ be a $\r
C^1$ (possibly, piecewise $\r C^1$) curve in $S^1_R$, given with
an orientation.
 Let us assume that its length $s(t_1,t_2)$ between the time
moments $t_1$ and
$t_2$ is  $>2\pi R\cdot \theta$. Then there exist time moments
$\tau_1$ and $\tau_2, \ t_1 < \tau_1 < \tau_2 < t_2$, such that
the following is true:}
\ev
\item{(i)}  {\sl The points $\omega (\tau_1)$ and $\omega (\tau_2)$
coincide or are antipodal points on $S^1_R$.}

\item{(ii)} {\sl  The velocity vectors $ v(\tau_i)= \D {d \omega \over dt}
(\tau_i), \ i=1,2$ satisfies:

\hskip0,7cm $\bullet$ \ $ \Vert v(\tau_i) \Vert \ge \D{1 \over 4} \D
{s(t_1,t_2)\over t_2-t_1}$ in the case of a closed curve,

\hskip0,7cm $\bullet$ \ $ \Vert v(\tau_i) \Vert \ge
\D { \theta- 4\over
4\theta } \cdot \D
{s(t_1,t_2)\D\over t_2-t_1}$ in the case of a non-closed curve,

\hskip0,7cm $\bullet$ \ the
directions of $v(\tau_1)$ and $v(\tau_2)$ are opposite
on the tangent line of $S^1_R$ at $ \omega(\tau_1)$ and
$ \omega(\tau_2)$ (in case $\omega(\tau_1)$ and
$ \omega(\tau_2)$ are not equal but antipodal, $v(\tau_2)$ is identified
with its projection on the tangent line at $\omega(\tau_1)$). }
\ev
{\bf Remark. }
Since for any $\theta>4$ we always have $\D{\theta-4\over
4\theta}\le {1\over 4}$,  in Lemma $3.6.$(ii) we have in both case
(the curve being closed or not): $ \Vert v(\tau_i)\Vert \ge \D {
\theta - 4\over 4\theta}\cdot \D {s(t_1,t_2)\D\over t_2-t_1}$, as
required by Lemma $3.5$. \ev
{\bf Proof of Lemma 3.6.}
The curve $\omega $ can be considered as ``lying over" the circle
$S^1_R$. Let  $s \in S^1_R$, $t \in {\omega }^{-1}(s)$ and let us
  denote $\alpha = \D{4\over \theta -4}$ and $T=t_2-t_1$
 the time interval.

First of all, we can always reduce the situation to the case of a
closed curve $\omega $. Indeed, if $\omega $ is not closed, but is
contained in $S^1_R$, we just ``close up" this curve with a
circular curve joining the endpoints, passed with the velocity $\D
{2\pi R\over \alpha T}$. We denote the new length of our curve
after the ``closing up" process by $\ti L $. We observe that the
new time interval $\ti T$ satisfies: $\ti T\le (1+\alpha)T$.
Assuming that Lemma 3.6 is true {\it for closed curves}, for the
new closed curve we therefore find some points with the velocity
at least:
 $$\D {\tilde L\over 4\tilde T}\ge {s(t_1,t_2)\over 4\ti T}\ge
{s(t_1,t_2)\over 4(1+\alpha) T}=\D { \theta - 4\over 4\theta}\cdot
\D {s(t_1,t_2)\D\over t_2-t_1}.$$ On the other hand, by
assumptions we have $${s(t_1,t_2)\over 4(1+\alpha) T}>{2\pi R\cdot
\theta\over 4(1+\alpha) T} = {2\pi R\over \alpha T}.$$ Hence, by
the above choice of the velocity on the added interval, the points
found have to belong to the original curve, and not to the added
interval.\ev

So from now on we assume that $\omega (t_1)=\omega (t_2)$ and
since we are allowed to deal with antipodal points and compare the
velocity vectors at these points, we can assume in what follows
that $\omega$ is homotopically trivial in $S^1_R$. So it can be
considered as ``lying over $\d  R$". Remind that the total length
of $\omega $ over a subset $A \subset \d  R$ is equal to the
integral $\int_A N(s) ds$ of the number $N(s)$ of the points of
$\omega $ over $A$.
 The same is true for  the ``positive" and the
``negative length" of $\omega$, since $\omega$ is homotopically
trivial in $S^1_R$ and thus may be view as a closed curve in a
segment. On the other hand, since $\omega (t_1)=\omega (t_2)$,
$\omega $ covers exactly the same length both in the positive and
in the negative directions. The same is true over any subset $A
\subset \d  R$. Therefore, over each subset $A \subset \d  R$ the
``positive length" of $\omega $ is equal to its ``negative
length".

\ev

Now let us assume that the statement of Lemma $3.6$ is not true for
the curve $\omega $. Denote its
length $s(t_1,t_2)$ by $L$. We thus assume there is no point in $w$
over which the velocities in the opposite directions are both
$\ge L/4T $. Consequently,
  over each point $s \in S^1_R$ either  all the positive
velocities of $\omega $ are   $<\D L/ 4T$, or all the negative
velocities of $\omega $ are  $<\D L/ 4T$, or both.

Denote by $A_1$ (respectively, $A_2, \ A_3$) the sets of points in
$S^1_R$ where the first (respectively, the second or the third)
alternative holds. Let us take those of the sets $A_1, \ A_2$ over
which the total length of $\omega $ is larger, say, $A_1$, and let
$A=A_1 \cup A_3$. The total length of $\omega $ over $A$ is $\ge
L/ 2$.

\ev

By the remark above, the ``positive length" of $\omega $ over $A$
is equal to its ``negative length", and hence each is $\ge L/ 4$.
But by the construction, at each positive point of $\omega $ over
$A$ the velocity is  $<L/ 4T$. Therefore, the total time required
for $\omega $ to cover $A$ in the positive direction is  $\D >{L/
4 \over  L/4T} = T$. This contradiction proves Lemma $3.6$, thus
it also proves Lemma $3.5$ and finishes the proof of Theorem
$3.4$.\f \ev The Euclidean version of Lemma $3.5$ is obtained by
using the Euclidean integral-geometric Cauchy-Crofton formula ([Fe
1] 5.11, [Fe 2] 2.10.15).

 Let us define the constant $C_n$ by $C_n=c_nV_n$, where $c_n$ is
the constant in the Cauchy-Crofton formula for
curves in $\d R^n$, and $V_n$
is the volume of the unit sphere in $\d R^n$. We have explicitly
$c_n=\D{ \Gamma ({n+1 \over 2}) \Gamma ({1\over 2}) /
\Gamma ({n \over 2})}$,
and $V_n=\D{2\Gamma^n(\D {1\over 2}) / \Gamma(\D {n\over 2})}$,
 where $\Gamma $ is the Euler
function. With these notations, let us state the Euclidean version
of Lemma $3.5$. \ev
{\bf Proposition  3.7. --- }
{\sl Let $\omega$ be a (piecewise) $\r C^1$ curve in a ball of
radius $R$ of $\d R^n$ and let, $\theta >8$ be a real number. If
the length $ s(t_1,t_2) $ of the curve $\omega$ between $t_1$ and
$t_2$ satisfies $s(t_1,t_2)>\theta\cdot C_n \cdot R$, then
 there exist a straight
line $\lambda$ and time moments $\tau_1$ and $\tau_2, \ t_1 <
\tau_1 < \tau_2 < t_2$, such that the following is true:}

\ev
\item{(i)} {\sl The orthogonal projections $\pi(\omega(\tau_1))$ and
$\pi(\omega (\tau_2))$ of the points $\omega(\tau_1)$ and
$\omega (\tau_2)$ on the line $\lambda$ coincide.}

\item{(ii)} {\sl
The norms of the orthogonal projections of the velocity vectors
$v(\tau_i)=\D {d\omega  \over dt} (\tau_i), \ i=1,2$ on the line $\lambda$
are   $>\D{\theta -8 \over 4\theta} \cdot  {s(t_1,t_2) \over t_2-t_1}$,
and the directions of these projections are opposite. }

\Ev
{\bf 3.3. Rotation of two Lipschitz trajectories}
\ev
 In this section we prove that the absolute rotation of any two
trajectories of a Lipschitz vector field in $\d R^3$ is bounded in
terms of the Lipschitz constant $K$ and the time interval.

\ev
{\bf Theorem 3.8. --- }
{\sl  Let $\omega_1$ and $\omega_2$ be trajectories,
on time intervals $T_1$ and $T_2$ respectively,
of a Lipschitz vector field $v$ defined in some open
subset $U$ of ${\d  R}^3$.
Then the mutual absolute rotation of $w_1$ and $w_2$ satisfies:
$$R_{abs}(\omega_1,\omega_2)\le  {K\over \pi}\cdot
\hbox{\tm min} (T_1,T_2) + {1\over 4\pi}\cdot K^2\cdot T_1\cdot T_2,$$
where $K$ is the Lipschitz constant of $v$. }
\ev

In fact, we shall prove a more accurate version of this theorem,
which  bounds the absolute rotation $R_{abs}(\omega_1,\omega_2)$
through the Lipschitz constant of the field and through the bound
of the absolute rotation of the trajectory $\omega_1$
(respectively $\omega_2$) around the {\sl points} of $\omega_2$
(respectively of $\omega_1$) (Theorem $3.9$). Then to get back to
Theorem $3.8$ we use the uniform bound on the rotation of
Lipschitz trajectories around points, provided by Theorem $3.4$.

\ev
{\bf Theorem 3.9. --- }
{\sl Let $\omega_1$ and $\omega_2$ be trajectories, on time
intervals $T_1$ and $T_2$ respectively, of a Lipschitz vector
field $v$ defined in some open subset $U$ of
${\d  R}^3$. With the following notations:
$$ R_1=\D \hbox{\tm max}_{p_2\in w_2}  R_{abs}(\omega_1,p_2), \
\  R_2=\D  \hbox{\tm max}_{p_1\in w_1}
 R_{abs}(\omega_2,p_1),$$
 the absolute rotation of $w_1$ and $w_2$ satisfies:
$$ R_{abs}(\omega_1,\omega_2)\le \D {1\over 4\pi }\cdot
  K \cdot R_1 \cdot  T_2
\ \hbox{ and } \
R_{abs}(\omega_1,\omega_2)\le \D {1\over 4\pi } \cdot K \cdot R_2\cdot  T_1.$$ }
\ev
{\bf Corollary 3.10. --- }
{\sl With the same notations and hypothesis as in Theorem $3.9$,
we have:
$$R_{abs}(\omega_1,\omega_2) \leq \D{K\over 4\pi} \cdot\min \{R_1\cdot
T_2, R_2 \cdot T_1 \}.$$ }
\ev
{\bf Proof of Theorem 3.9.}
We recall that for trajectories of the vector field $v$ the Gauss
integral takes the form:
$$ R_{abs}(\omega_1,\omega_2) =
{1\over 4\pi} \int_{T_1}
\int_{T_2}\ {\vert <v_1\wedge v_2,r_{12}>\vert  \over \Vert r_{12}\Vert^3} \
\  dt_1 dt_2 . $$

\noindent
Here $v_1,v_2$ are the velocity vectors of $\omega_1,\omega_2$,
respectively,
(i.e. the values of the vector field $v$ at the running points $p_1(t_1)$ and
$p_2(t_2)$  on $\omega_1$
and $\omega_2$), and $r_{12}=p_2-p_1$ is the vector
joining the running points $p_1$
and $p_2$. So we have:

$$4\pi \cdot R_{abs}(\omega_1,\omega_2) = \int_{T_1} \int_{T_2}\
{\vert<v_1\wedge v_2,r_{12}>\vert \over \Vert r_{12}\Vert^3} \ \
dt_1 dt_2  $$
$$=  \int_{t_2\in T_2} dt_2 \int_{t_2\in T_1}
{\vert<v_1(t_1)\wedge (v_2(t_2) - v_1
(t_1)),r_{12}(t_1,t_2)>\vert \over \Vert r_{12}\Vert^3}\ \ dt_1.
$$
  Indeed, the substraction of $v_1(t_1)$ from $v_2(t_2)$
does not change the triple product under the integral.  Now since
the vector field $v$ is Lipschitz, we have: $ \Vert
v_2(t_2)-v_1(t_1)\Vert\leq K\cdot \Vert r_{12}\Vert .$

\n
Hence for the triple product we obtain:
$$ \vert<v_1\wedge(v_2-v_1),r_{12}>\vert\leq\Vert
\tilde v_1\Vert\cdot K\cdot \Vert r_{12}
\Vert\cdot\Vert r_{12}\Vert , $$
where $\tilde v_1$ is the component of the vector $v_1$ orthogonal to the
vector $r_{12}$ (see Figure $6$).

\vskip6.3cm
\vglue0cm
\vskip0cm
\hskip0cm\includegraphics{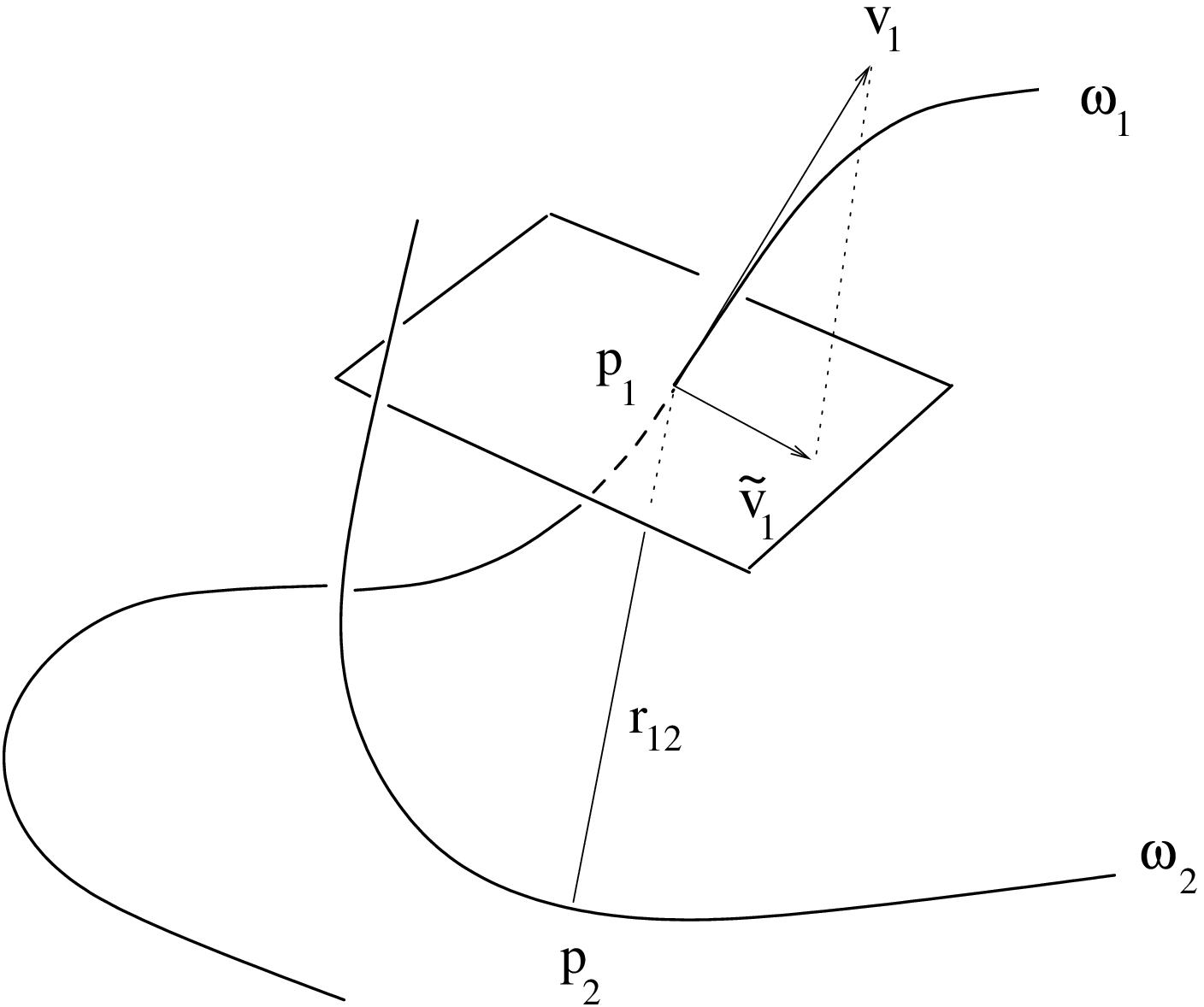}
\vglue0cm
\Ev
\centerline{\bf Figure 6}
\Ev

\n
 Therefore for the absolute rotation we get:
$$4\pi\cdot  R_{abs}(\omega_1,\omega_2) \leq K \cdot
\int_{t_2\in T_2}dt_2 \int_{t_1\in T_1}\
{\Vert \tilde v_1(t_1,t_2)\Vert \over \Vert r_{12}(t_1,t_2)\Vert}\ \
 dt_1 . $$
In the interior integral  the time $t_2$ and the corresponding
point $p_2(t_2)$ on the curve $\omega_2$ are fixed, while the
integration runs over the curve $\omega_1$, hence the integral
$\D\int_{t_1\in T_1} {\Vert \tilde v_1(t_1,t_2)\Vert \over \Vert
r_{12} (t_1,t_2)\Vert}\ \ dt_1$ is equal to the length of the
spherical projection $\sigma_1$ of the curve $\omega_1$ from the
point $p_2(t_2)$, i.e. to $2\pi\cdot R_{abs}(\omega_1,p_2)$. By
the assumptions, this rotation is uniformly in $p_2$ bounded by
$R_1$. Hence the interior integral $\D \int_{t_1\in T_1}{\Vert
\tilde v_1\Vert \over \Vert r_{12}\Vert}\ \ dt_1$ does not exceed
$R_1$, and finally:
$$4\pi R_{abs}(\omega_1,\omega_2) \leq K \cdot T_2 \cdot R_1.$$
The setting of Theorem $3.9$ is symmetric with respect to the
trajectories $\omega_1$ and $\omega_2$. So interchanging these
trajectories we get:
$$4\pi R_{abs}(\omega_1,\omega_2) \leq K \cdot T_1 \cdot R_2.$$
This completes the proof of Theorem 3.9.\f

\ev
{\bf Proof of Theorem 3.8.}
 We use Theorem $3.4$,
which states that a rotation in time $T$ of any trajectory of a
Lipschitz vector field $v$ in ${\d  R}^3$ around any point $p$
does not exceed:
$$ 4+ K\cdot T. $$
Applying Corollary $3.10$ we get:
$$ R_{abs}(\omega_1,\omega_2) \leq
\hbox{\tm min} ( {K\over \pi}\cdot T_1+
{1\over 4\pi}\cdot K^2\cdot T_1\cdot T_2,
{K\over \pi}\cdot T_2+
{1\over 4\pi}\cdot K^2\cdot T_1\cdot T_2,
),$$
and Theorem  $3.8$ is proved.\f
\ev {\bf Example of application: Logarithmic bound of the local
rotation at singular points for analytic vector  fields.} \ev

Let us consider the following situation. The vector field:
$${dv\over dt} = Lv+G(v)$$ is defined in a neighborhood of
$O\in {\d
R}^3$ and has a non-degenerate linear part $L$ with all the
eigenvalues $\ell_j, \ j=1,2,3,$ having  a negative real part:
$\Re (\ell_j) \leq \ell <0, \ j=1,2,3.$ In dynamical language, $v$ has
a non-degenerate sink at the origin (it is the case of
our field in $\d R^2$,  in the Remark following Proposition $3.1$).
It is easy to see that the
Lipschitz constant of $v$ in a neighborhood $U$ of the origin
tends to the norm of $L$ as $U$ shrinks to the origin. The
following theorem is an immediate corollary of the results of
Section $3$:
\ev
{\bf Theorem 3.10. --- }
{\sl For any two trajectories $\omega_1, \omega_2$ of the field
$v$ in a neighborhood of $O\in {\d  R}^3$, the absolute
rotation $R_{abs}(\omega_1, \omega_2)$ grows at most logarithmically
with
the distance to the origin. More accurately, the rotation
$R_{abs}(\omega_1, \omega_2, R, r)$ of the parts of
$\omega_1, \omega_2$ between
the spheres of the radii $R > r > 0$ satisfies:
$$ R_{abs}(\omega_1,
\omega_2, R, r) \leq C \Vert L \Vert {{\log^2 ({R/ r})}\over \ell}.$$
}
\ev
{\bf Proof.}
This is a direct consequence of Theorem $3.8$, since the
Lipschitz constant of $v$ in a neighborhood $U$ of the origin
tends to $\Vert L \Vert$ as $U$ shrinks to the origin, while the
time interval for both the trajectories between the spheres of the
radii $R > r > 0$ is of order $\D {\log ({R/ r})}\over \D \ell$.
\f
\ev
{\bf Remark.}
Of course, one can easily show that the bound of Theorem $3.10$ is
sharp: consider a linear vector field:
$$ {dv\over dt} = L v $$
in a neighborhood of $O\in {\d  R}^3$, with a non-degenerate
linear part $L$ having all its eigenvalues $\ell_j, \ j=1,2,3,$ with
 negative real part: $Re (\ell_j) \leq \ell <0, \ j=1,2,3.$
Assume in addition that $\ell_1 \in \d  R$, while $\ell_2$ and
$\ell_3$ are conjugate: $\ell_{2,3}= \alpha \pm i \beta$. Then the
solutions are $x_1=C_1 \cdot \exp (\ell_1 t), \ x_2=C_2 \cdot
\exp( \alpha t )\cdot\sin
(\beta t), \ \ x_3=C_3 \cdot\exp (\alpha t )\cdot
\cos (\beta t),$ and the
trajectories rotate around the $x_1$-axis and one around another
exactly as prescribed by the upper bound given in Theorem $3.8$.

\Ev
\Ev
\centerline{\bf References }
\Ev

\item{[Ar-Kh]}
 Arnold, V. I.; Khesin, B. A. Topological methods
in hydrodynamics. {\sl Applied Mathematical
Sciences,} {\bf 125}
Springer-Verlag, New York, 1998. xvi+374 pp.

\vskip1mm
\item{[Bl-Mo-Ro]} Blais, F.; Moussu, R.; Rolin, J.-P.
Non-oscillating integral curves and o-minimal structures.
{\sl Analyzable functions and applications, 103--112,
Contemp. Math.,} {\bf 373},
Amer. Math. Soc., Providence, RI, 2005.

\vskip1mm \n
\item{[Ca-Mo-Ro]}
Cano, F.; Moussu, R.; Rolin, J.-P.
Non-oscillating integral curves and va\-lua\-tions.
{\it J. Reine Angew. Math.} {\bf 582}, 2005, 107--141.

\vskip1mm \n
\item{[Ca-Mo-Sa 1]}
 Cano, F.; Moussu, R.; Sanz, F.
Oscillation, spiralement, tourbillonnement.
{\sl  Comment. Math. Helv.} {\bf 75}, 2000, no. 2,
284--318.

\vskip1mm
\item{[Ca-Mo-Sa 2]} Cano, F.; Moussu, R.; Sanz, F.
 Pinceau int\'egral enlac\'e.
{\sl C. R. Math. Acad. Sci. Paris} {\bf 334},
 2002, no. 10, 855--858.
 \vskip1mm

\item{[Ca-Mo-Sa 3]} Cano, F.; Moussu, R.; Sanz, F.
Pinceaux de courbes int\'egrales d'un champ de vecteurs analytique.
{\sl  Analyse complexe, syst\`emes dynamiques, sommabilit\'e des s\'eries
divergentes et th\'eories galoisiennes. II. Ast\'erisque}
 {\bf 297},  2004, 1--34.

 \vskip1mm
 \item{[Co]}
Coste, M.
An introduction to o-minimal geometry.
{\sl Publ. Dipartimento Di Matematica Dell'Universi\`a di Pisa.}
(2000)

 \vskip1mm
 \item{[Dr]}
van den Dries, L.
Tame topology and o-minimal structures.
{\sl London Mathematical Society Lecture Note Series}, {\bf 248}.
Cambridge University Press, Cambridge, (1998), x+180 pp.

 \vskip1mm
\item{[Dr-Mi]}
  van den Dries, L.; Miller, C. Geometric categories
and o-minimal structures. {\sl Duke Math.
 J.} {\bf 84} (1996), no. 2, 497--540.

\vskip1mm \n
\item{[Du-Fo-No]}
Dubrovin, B. A.; Fomenko, A. T.; Novikov, S. P. Modern
geometry - methods and applications. Part II. The geometry and
topology of manifolds. Translated from the Russian by Robert
G. Burns. {\sl Graduate
Texts in Mathematics} {\bf 104}, Springer-Verlag, New York, 1985.
xv+430 pp.

  \vskip1mm
\item{[Fe 1]} H. Federer, The $(\Phi,k)$ rectifiable subsets
of n space. {\it Trans. Amer. Math. Soc.} {\bf 62},
 (1947), 114-192.

 \vskip1mm
\item{[Fe 2]} H. Federer, Geometric measure theory. {\it
Grundlehren Math. Wiss.} {\bf 153}, Springer Verlag,
(1969).

 \vskip1mm
 \item{[Gr-Yo]} Grigoriev, A; Yomdin, Y.
Rotation rate of a trajectory of an algebraic vector field
around an algebraic curve.  {\sl
Qual. Theory Dyn. Syst.} {\bf 2}, 2001, no. 1, 61--66.

  \vskip1mm
\item{[Ho 1]} Hovanskii, A. G.
Fewnomials.
{\sl Translations of Mathematical Monographs} {\bf 88}.
{\sl American Mathematical Society, Providence, RI,} 1991.
viii+139 pp.

  \vskip1mm
\item{[Ho 2]} Hovanskii, A. G.
Fewnomials and Pfaff manifolds.
{\sl
Proceedings of the International Congress of Mathematicians,
 Vol. 1, 2 (Warsaw, 1983)},
549–564, PWN, Warsaw, 1984.
  \vskip1mm

\item{[Ho-Ya]}
 Hovanskii, A.; Yakovenko, S. Generalized Rolle theorem in
 $\d R^ n$ and $\d C$. {\sl J. Dynam. Control Systems } {\bf 2},
 1996, no. 1, 103--123.
\vskip1mm

\item{[Ku]} Kurdyka, K.
On gradients of functions definable in o-minimal
structures. {\sl Ann. Inst. Fourier} {\bf 48}
 (1998), no. 3, 769--783.

\item{[Ku-Mo]} Kurdyka, K; Mostowski, T.
On the gradient conjecture of Ren\'e Thom.
 {\sl  Preprint Universit\'e de Savoie,} (1996)

  \vskip1mm

\item{[Ku-Mo-Pa]} Kurdyka, K; Mostowski, T.; Parusinski A.
 Proof of the gradient conjecture of R. Thom. {\sl Ann.
of Math.} (2) {\bf 152} (2000), no. 3, 763--792.

  \vskip1mm

\item{[Li-Mo-Ro]} Lion, J.-M.; Moussu, R.; Sanz, F.
 Champs de vecteurs analytiques et champs de gradients.
 {\sl  Ergodic Theory Dynam. Systems} {\bf  22},
 (2002), no. 2, 525--534.
  \vskip1mm

\item{[No-Ya 1]}
Novikov, D.; Yakovenko, S. Integral curvatures, oscillation and
rotation of spatial curves around
affine subspaces.
{\sl J. Dynam. Control Systems} {\bf 2}, 1996, 157--191.
    \vskip1mm

 \item{[No-Ya 2]}   Novikov, D.; Yakovenko, S.
Meandering of trajectories of polynomial vector fields in the
affine $n$-space. {\sl Proceedings of the Symposium on Planar
 Vector Fields (Lleida, 1996). Publ. Mat.} {\bf 41},
 1997, no. 1, 223--242.
       \vskip1mm

  \item{[Sh]}
 Shiota, M. Geometry of subanalytic and semialgebraic sets.
{\sl Progress in Mathematics},
{\bf 150} Birkhäuser Boston, Inc., Boston, MA, (1997),
 xii+431 pp.

   \vskip1mm
  \item{[Wi]} Wilkie, A. J.
A theorem of the complement and some new o-minimal structures
{\sl Selecta Math. (N.S.)} {\bf 5}, 1999, no. 4, 397–421.
     \vskip1mm

\item{[Ya]} Yakovenko, S. On functions and curves defined by ordinary
differential equations. The Arnoldfest (Toronto, ON, 1997),
497--525, {\sl Fields Inst. Commun.} {\bf 24}
, {\sl Amer. Math. Soc., Providence,
RI,} 1999.
   \vskip1mm

\end